\magnification=\magstep1
\input amstex
\documentstyle{amsppt}
\catcode`\@=11 \loadmathfont{rsfs}
\def\mycal{\mathfont@\rsfs}
\csname rsfs \endcsname \catcode`\@=\active  

\vsize=7.5in

\topmatter 
\title 
Compressible subalgebras in II$_1$ factors \\ 
$\text{\it To Huzihiro Araki, in memoriam}$
\endtitle
\author  Sorin Popa   \endauthor

\rightheadtext{Compressible subalgebras}

\affil University of California Los Angeles, popa\@math.ucla.edu \endaffil

\address Math.Dept., UCLA, LA, CA 90095-155505\endaddress
\email popa\@math.ucla.edu\endemail

\thanks Supported in part by the Takesaki Endowed Chair in Operator Algebra at UCLA (4/23/26) \endthanks

\abstract Given a II$_1$ factor $M$, a W$^*$-subalgebra $Q\subset M$ is {\it compressible} 
if  for any $\varepsilon>0$ there exists a finite set of unitary elements $\Cal U_0\subset \Cal U(M)=\Cal U(M\otimes 1)$ such 
that $\| \frac{1}{|\Cal U_0|}\sum_{u\in \Cal U_0} uxu^* -E_{1\otimes \Bbb M_K(\Bbb C)}(x)\|\leq \varepsilon$, 
$\forall K\geq 1$, $\forall x\in (Q\otimes \Bbb M_K(\Bbb C))_1$.  Any W$^*$-subalgebra $Q$ in a II$_1$ factor $M$ 
which admits a diffuse W$^*$-algebra $Q_0\subset M$  that's free independent to $Q$, is compressible in $M$. We prove that if 
$Q\subset M$ is compressible, then $_NL^2M_Q$ contains a copy of the coarse $N-Q$ bimodule for any AFD subalgebra $N\subset M$. 
We use this result to provide examples of inclusions of factors $M\subset \Cal M$ that are ergodic but not AFD-ergodic, even after stabilizing by $\Cal B(\ell^2\Bbb N)$. 

\endabstract

\endtopmatter

\document

\heading 0.  Introduction  \endheading





We consider in this paper a new property for an inclusion of von Neumann algebras (hereafter called a W$^*$-inclusion) $Q\subset M$:  
we say that $Q\subset M$ is {\it compressible} if  
the averaging action of the unitary group $\Cal U(M)$ on matrices over $M$ pushes 
all matrices of norm $\leq 1$ with entries in $Q$  to scalar matrices, uniformly in operator norm. 
That is, for any $\varepsilon>0$ there exists a finite set of unitary elements $\Cal U_0\subset \Cal U(M)$ such 
that $\| \frac{1}{|\Cal U_0|}\sum_{u\in \Cal U_0} uxu^* -E_{1\otimes \Bbb M_K(\Bbb C)}(x)\|\leq \varepsilon$, 
$\forall K\geq 1$, $\forall x\in (Q\otimes \Bbb M_K(\Bbb C))_1$.  

We use a result in [PV14] to show that if $M$ is a II$_1$ factor and $Q\subset M$ is a diffuse W$^*$-subalgebra  which admits 
some $u\in \Cal U(M)$ 
with $\tau(u)=0$ and $\{u, u^*\}$ free independent to $Q$, then $Q\subset M$ is compressible. Since  the existence of such $u$ 
implies $uQu^*$ free independent to $Q$, this  condition is in fact equivalent to requiring the existence of a diffuse W$^*$-subalgebra $Q_0\subset M$ that's free independent to $Q$. 

Our main result shows that compressible subalgebras are in some sense ``AFD-repellent''. Namely,  if $Q\subset M$ is compressible 
then, given any  approximately finite dimensional W$^*$-subalgebra (abbreviated {\it AFD-subalgebra}) $N\subset M$, the $N-Q$ bimodule 
$L^2M$ contains a copy of the coarse $N-Q$ bimodule $L^2N\overline{\otimes} L^2Q$. If in addition $N$ 
is quasi-regular in $M$, then $N$ is coarse to $Q$, i.e, $_NL^2M_Q\subset (L^2N\overline{\otimes} L^2Q)^{\oplus\infty}$.

\proclaim{0.1. Theorem} Let $Q\subset M$ be a tracial W$^*$-inclusion. If $Q\subset M$ is compressible, then given any AFD-subalgebra  $N \subset M$,  the Hilbert bimodule $_NL^2M_Q$ has a non-zero coarse part. 
If in addition $N$ is quasi-regular in $M$, then  $N$ is coarse to $Q$. 
\endproclaim


This easily implies that a diffuse tracial quasi-regular W$^*$-inclusion, $Q\subset M$, cannot be compressible   
and that if a tracial W$^*$-algebra $M$ contains a diffuse compressible W$^*$-subalgebra then it has no AFD direct summand.  In particular,  
the hyperfinite II$_1$ factor has no compressible diffuse subalgebras,  
and if  a group $G$ has an infinite subgroup $H$ such that $LH\subset LG$ is compressible  then $G$ is not amenable. 

The theorem shows  that  compressibility is a strong obstruction to tightness 
in II$_1$ factors. Recall from (page 993 in [P19a]) that a pair of 
factors $N_0, N_1$ in a II$_1$ factor $M$ is   {\it tight} if $N_0\vee N_1^{op}=\Cal B(L^2M)$, or equivalently if 
the Hilbert bimodule $_{N_0}L^2M_{N_1}$ is irreducible. Also, given a subfactor $Q\subset M$, 
a subfactor $N\subset M$ with the property that $_N L^2M_Q$ is irreducible is called a {\it tight complement} 
of $Q$. As shown in (Section 7 of [P19a]; cf. also [P19b]), 
of particular interest are the tight pairs and tight complements involving 
the hyperfinite II$_1$ factor.  Since an irreducible bimodule cannot contain a coarse bimodule, the above result implies 
that given any embedding of a free product 
$Q*Q_0$ into a II$_1$ factor $M$, with $Q$ a factor and $Q_0\neq \Bbb C$, the subfactor $Q\subset M$ cannot have a hyperfinite tight complement.

This is also related to Kadison's problem, asking whether an ergodic (or irreducible) W$^*$-inclusion of factors $\Cal N\subset \Cal M$ 
is necessarily MASA-ergodic, i.e., $\Cal N$ contains a maximal abelian $^*$-subalgebra (abbreviated MASA) of $\Cal M$. 
By [P81], this problem has an affirmative answer in case $\Cal N, \Cal M$ are II$_1$, and more generally when  $\Cal N$ is type II 
and there exists a normal conditional expectation of $\Cal M$ onto $\Cal N$ (see [M23] for a general 
result along these lines). But it has been shown in [GP96] that it fails for 
any basic construction inclusion of factors $M\subset M_1=\langle M, e_Q\rangle$, with $M$ a free group factor and $Q\subset M$ an ergodic 
embedding of the hyperfinite II$_1$ factor.  
The result in [GP96] showed in fact more, namely that a free group factor $M$ cannot be decomposed 
as $M=\overline{\text{\rm span}NQ}^w$ for some AFD-subalgebras $N, Q\subset M$. In particular, if we consider $M\subset M_1= \langle M, e_Q\rangle$ 
for some hyperfinite II$_1$ subfactor $Q\subset M$, then this inclusion is not AFD-ergodic, i.e., there exists no 
AFD-subalgebra $N\subset M$ such that $N'\cap M_1=\Cal Z(N)$. In particular, it is not $R$-ergodic, i.e., 
there exist no hyperfinite subfactors $N\subset M$ such that $N'\cap M_1=\Bbb C$.

We show that even after stabilizing by $\Cal B(\ell^2\Bbb N)$, a basic construction inclusion   $M\subset \langle M, e_Q\rangle$ 
arising from a compressible subfactor $Q\subset M$ cannot be AFD-ergodic, in particular it is not MASA-ergo (thus 
failing Kadison's problem).

\proclaim{0.2. Theorem} Let  $Q\subset M$ be an ergodic inclusion of $\text{\rm II}_1$ factors and denote 
$M\subset M_1=\langle M, e_Q \rangle$ the associated basic construction W$^*$-inclusion. If $Q\subset M$ is compressible 
$($e.g., if $M\supset Q*Q_0$ for some diffuse W$^*$-algebra $Q_0)$, 
then $M\subset M_1$ is ergodic but not AFD-ergodic. In particular, it is not MASA-ergodic, nor $R$-ergodic. 
Moreover, $M\overline{\otimes} \Cal B(\Cal H_0)\subset M_1\overline{\otimes} \Cal B(\Cal H_0)$ 
is not AFD-ergodic  for any Hilbert space  $\Cal H_0$. 
\endproclaim

The paper is organized as follows. In Section 1 we introduce the compressibility property, prove its basic properties and give examples. 
In Section 2 we prove Theorem 0.1 (see Corollary 2.9) and Theorem 0.2 (see Corollary 2.10). They will be derived from more technical 
general results stated as Theorems 2.3,  2.5, whose proofs are based on repeated usage of Kaplansky's density theorem. 

Besides its relevance to tightness problems in II$_1$ factors, our  study  of compressibility 
for W$^*$-inclusions was motivated by the free complementation 
problem in free group factors, and its weaker version 
asking whether any AFD-subalgebra  $Q$ of $L\Bbb F_n$ admits a 
diffuse subalgebra  $Q_0\subset L\Bbb F_n$ that's free to $Q$. An initial motivation came also from Connes embedding problem. 
We comment on this  and on other aspects in Section 3, where we also 
formulate some open problems.   

\vskip.05in 
{\it Acknowledgement}. This paper is dedicated to the memory of Huzihiro Araki, whose pioneering work in operator 
algebras and mathematical physics had a deep, lasting impact on these subjects. The Araki-Woods AFD factors and their 
free version, conceived by Dima Shlyakhtenko in the 1990s, continue to provide the most important  class of examples 
of type III factors, a fundamental framework for ``non-commutative analysis''. His 1971 work with 
Mi-Soo Bae Smith and Larry Smith on the homotopy of the unitary group $\Cal U(M)$ of a von Neumann algebra $M$ endowed 
with the operator norm directly inspired Masamichi Takesaki and myself to study in the 1990s 
the contractibility of $\Cal U(M)$ when $M$ is a II$_1$ factor 
and $\Cal U(M)$ is endowed with the Hilbert norm $\| \ \|_2$ given by the trace. This problem has recently seen a 
spectacular revival through work of Narutaka Ozawa early 2025, leading to the final striking resolution by David Jekel in the Summer of 2025. 
Most recently, the Araki-Smith-Smith result that $\pi_1(\Cal U(M))=\Bbb Z$ for II$_1$ factors $M$ played 
a crucial role in the brilliant new proof by David Gao and Srivatsav Kunnawalkam Elayavalli of the famous Pimsner-Voiculescu 
1981 result that the reduced C$^*$-algebras of free groups of different rank are non-isomorphic. Araki's formidable 
energy has been as inspiring for me and all mathematicians who had the chance 
to interact with him, as was his mathematical work. He will be thoroughly missed. 

\heading 1.  Compressibility for subalgabras  \endheading 

For basic facts on II$_1$ factors  we send the reader to [AP17] and for general operator algebras to [S72], [T79], [BrO08]. 
This paper is much related to our previous work in  
[P18], [P19a] and [P19b], from which we'll use terminology and style of notations. 

As usual, inclusions of von Neumann algebras will be unital, and will be often called 
W$^*$-inclusions, W$^*$-subalgebras, etc. 

Given a normed space $(\Cal X, \| \ \|)$, we denote the unit ball of $\Cal X$ by $(\Cal X)_1$.  

\vskip.05in

\noindent 
{\bf 1.1. Definitions.} Let $\Cal M$ be a W$^*$-algebra with W$^*$-subalgebras $Q, B \subset \Cal M$. In what follows we identify 
$\Cal M$ with the subalgebra $\Cal M \otimes 1$ in $\Cal M \otimes \Bbb M_K(\Bbb C)$, $K\geq 1$. 

$1^\circ$ If for some $\Cal U_0\subset \Cal U(\Cal M)$  finite and $\varepsilon >0$ we have that  
for any $K\geq 1$ and any $x\in (Q \otimes \Bbb M_K(\Bbb C))_1$,  
there exists $x_0 \in (B \otimes \Bbb M_K(\Bbb C))_1$ such that 
$\|\frac{1}{|\Cal U_0|} \sum_{u\in \Cal U_0} uxu^* -x_0\|\leq \varepsilon$, then we say that $Q\subset \Cal M$ is $(\Cal U_0, \varepsilon)$-{\it compressible } {\it relative to} $B$. 

$2^\circ$  We say that $Q\subset \Cal M$  
{\it is compressible} 
{\it relative to $B$}  if for any $\varepsilon >0$ there exists a finite set of unitaries $\Cal U_\varepsilon \subset \Cal U(\Cal M)$ such that 
$Q$ is $(\Cal U_\varepsilon, \varepsilon)$-compressible relative to $B$. 
In case this condition holds true for $B=\Bbb C1$, then we simply say that 
$Q\subset \Cal M$ is {\it compressible}. Also, if all $\Cal U_\varepsilon$ can be chosen in some $\Cal U\subset \Cal U(\Cal M)$, 
then we say that $\Cal U$ is a {\it compressing set of unitaries} for $Q\subset \Cal M$. 
  
 \vskip.05in 
 
 The following lemma provides some alternative ways of describing compressibility for subalgebras. 

\proclaim{1.2. Lemma} Let $Q, B\subset \Cal M$ be W$^*$-algebras, $\Cal U_0\subset \Cal U(\Cal M)$ a finite set and $\varepsilon >0$. 
The following conditions are equivalent: 
\vskip.05in
$(1)$ $Q\subset \Cal M$ is $(\Cal U_0, \varepsilon)$-compressible relative to $B$. 
\vskip.05in 

$(2)$ For any AFD W$^*$-algebra $\Cal B$, the inclusion $Q\overline{\otimes} \Cal B \subset \Cal M \overline{\otimes} \Cal B$ 
is $(\Cal U_0, \varepsilon)$-compressi- ble relative to $B\overline{\otimes} \Cal B$, 

$(3)$ For any AFD W$^*$-algebra $\Cal B$ and any $x\in (Q\overline{\otimes} \Cal B)_1$ there exists $x_0\in 
(B \otimes \Cal B)_1$ such that 
$$
\|\frac{1}{|\Cal U_0|}\sum_{u\in \Cal U_0} u x u^* -x_0\| \leq \varepsilon.    
$$

$(4)$ There exists a W$^*$-algebra $\Cal B$ that contains matrix subfactors of arbitrarily large dimension 
$($not necessarily with same unit as $\Cal B)$, such that for any $x\in (Q\overline{\otimes} \Cal B)_1$, there exists $x_0\in 
(B \otimes \Cal B)_1$ such that 
$$
\|\frac{1}{|\Cal U_0|}\sum_{u\in \Cal U_0} u x u^* -x_0\| \leq \varepsilon .   
$$
\endproclaim
\noindent
{\it Proof}. The implications $(2) \Rightarrow (3) \Rightarrow (4) \Rightarrow (1)$ are obvious. To prove $(1) \Rightarrow (2)$, 
we clearly only need to show that $(1)$ implies $(2)$ holds true for $\Cal B=\Cal B(\Cal H)$. But if $p_i\in \Cal B(\Cal H)$ 
is an increasing net of finite dimensional projections such that $p_i \nearrow 1$, then 
property $(1)$ implies that for any $x\in (Q \otimes p_i\Cal B(\Cal H)p_i)_1$ there exists $x_0 \in (B\otimes p_i\Cal B(\Cal H)p_i)_1$ 
such that 
$$
\|\frac{1}{|\Cal U_0|}\sum_{u\in \Cal U_0} u x u^* -x_0\| \leq \varepsilon.    
$$
Thus, we also have this property for all $x\in Q \otimes \cup_i p_i\Cal B(\Cal H)p_i$. 
If now $x\in (Q\overline{\otimes} \Cal B(\Cal H))_1$ is an arbitrary element and $x_j \in \cup_i Q\otimes p_i\Cal B(\Cal H)p_i$ is a net of elements 
with $\|x_j\|\leq 1$ and $x_j$ converging to $x$ in the $wo$-topology, then one gets elements $x_j^0\in (\cup_i B \otimes p_i\Cal B(\Cal H)p_i)_1$ 
such that 
$$
\|\frac{1}{\Cal U_0}\sum_{u\in \Cal U_0} u x_j u^* -x^0_j\| \leq \varepsilon   \tag 1.2.1 
$$
Taking the $wo$-limit in $(1.2.1)$ 
of the (bounded) net $\{x_j\}_j \in \cup_i Q\otimes p_i\Cal B(\Cal H)p_i \subset Q\overline{\otimes} \Cal B(\Cal H)$ 
and  respectively a $wo$-limit $x_0$ of the net 
$\{x_j^0\}_j \subset (\cup_i B \otimes p_i\Cal B(\Cal H)p_i)_1 \subset (B\overline{\otimes} \Cal B(\Cal H))_1$, 
and using the upper semi-continuity of the operator norm with respect to $wo$-topology, we thus get 
$$
\|\frac{1}{|\Cal U_0|}\sum_{u\in \Cal U_0} u x u^* -x_0\| \leq  \varepsilon.   
$$
\hfill $\square$

\proclaim{1.3. Lemma} $1^\circ$ Let $\Cal M$ be a W$^*$-algebra with W$^*$-subalgebras $Q, B\subset \Cal M$. Assume $Q\subset \Cal M$ 
is compressible relative to $B$. Then $Q_0$ is compressible relative to $B$ for any 
W$^*$-subalgebra $Q_0\subset Q$. If $\Cal M \subset \tilde{\Cal M}$ 
is an embedding of $\Cal M$ in a larger W$^*$-algebra, then $Q\subset \tilde{\Cal M}$ is compressible relative to $B$. 
Also, $Qp\subset \Cal Mp$ is compressible relative to $Bp$ for any $p\in \Cal Z(\Cal M)$. 

$2^\circ$ If $Q$ is a matricial factor, $Q\simeq \Bbb M_n(\Bbb C)$, then any W$^*$-embedding $Q\subset \Cal M$  is compressible and more generally 
any W$^*$-subalgebra $Q_0\subset Q$ is compressible in $\Cal M$.  

$3^\circ$ If $Q$ is diffuse and $Q\subset \Cal M$ is  compressible, then $\Cal M$ does not have any 
type $\text{\rm I}_{fin}$ direct summand. 

$4^\circ$ Let $Q\subset \Cal M$ be an inclusion of  von Neumann algebras. 
Assume $\{p_j\}_j\subset \Cal P(Q)$ is a finite 
partition of $1$ so that $\{p_j\}_j$ is contained in a matrix subfactor of $Q$ and such that  $p_jQp_j\subset p_j \Cal M p_j$ is compressible, $\forall j$. Then $Q \subset \Cal M$ is compressible.

$5^\circ$ Let $B\subset Q \subset \Cal M$ be inclusions of von Neumann algebras. If $Q\subset \Cal M$ is compressible relative to $B$ and $B \subset \Cal M$ 
is compressible, then $Q \subset \Cal M$ is compressible. 
\endproclaim
\noindent
{\it Proof}. $1^\circ$  is trivial from the definitions. 

2$^\circ$ If $\Cal U_0\subset \Cal U(Q)$ denotes the set of unitary matrices that have only $\pm 1$ and $0$ as entries, then $\Cal U_0$ is a finite subgroup 
of $\Cal U(Q)$ with $|\Cal U_0|=2^n n!$ and given any $x\in Q$ we have $\frac{1}{|\Cal U_0|} \sum_{u\in \Cal U_0} uxu^*=tr (x)1$, where $tr$ denotes 
here the normalized trace on $Q$. Thus, for any $K\geq 1$ and any $x\in Q\otimes \Bbb M_K(\Bbb C)$ we have  
$$
\frac{1}{|\Cal U_0|} \sum_{u\in \Cal U_0} uxu^*=E_{1\otimes \Bbb M_K(\Bbb C)}(x)
$$

$3^\circ$ Assume that $Q$ is diffuse and $Q\subset \Cal M$ is compressible. By cutting off with a central projection 
(see 1$^\circ$) we may assume $\Cal M$ of type I$_n$. Thus $\Cal M = \Cal Z \otimes \Bbb M_n(\Bbb C)$, with $\Cal Z$ abelian diffuse, 
which we identify here with the center of $\Cal M$. Since the averaging with unitaries preserves the central trace $ctr$ on $\Cal M$, 
$Q\subset \Cal M$ compressible implies that the restriction $ctr_{|Q}$ gives a normal state $\varphi$ on $Q$. But a normal state on a diffuse abelian 
W$^*$-algebra  has the property that there exist projections $q\in Q$ with $\varphi(q)$ arbitrary small, while on $\Cal M$ we have $ctr(q)\geq \frac{1}{n}q$ 
for any projection $q\in \Cal M$, contradiction. 

$4^\circ$ Let $\varepsilon > 0$. Since each $p_jQp_j\subset p_j \Cal M p_j$ is compressible, one can find finite subsets $\Cal U_j\subset \Cal U(p_j\Cal M p_j)$ 
such that the averaging by $\oplus \Cal U_j$ 
of any $x_j$ in the unit ball of $(p_j Q p_j)\otimes \Bbb M_K(\Bbb C)$ is less than $\varepsilon \tau(p_j)/2$-close to an element $y_j$ 
in the unit ball of $p_j 1 \otimes  \Bbb M_K(\Bbb C)$. By embedding $Q_0=\sum_j \Bbb C p_j$ in a matrix factor and applying $2^\circ$, 
one can thus get an averaging that takes $\sum_j y_j$ close to within $\varepsilon/2$ to a scalar multiple of $1$. Composing with the previous averager, 
one takes any $\sum_j x_j$ to within $\varepsilon$-close an element in $1\otimes \Bbb M_K(\Bbb C)$. For an arbitrary $x\in (Q \otimes \Bbb M_K(\Bbb C))_1$, 
one first take the averaging of $x$ by the set of unitaries of the form $\sum_j \alpha_j p_j$, where $\alpha_j=\pm 1$, which takes $x$ to 
$\sum_j x_j$, where $x_j=p_jxp_j$. 

$5^\circ$ This is similar to the proof of 4$^\circ$ above and we leave the details as an exercise.   

\hfill $\square$

\proclaim{1.4. Proposition}  Let $Q, Q_0$ be  tracial von Neumann algebras with a common von Neumann subalgebra $B\subset Q, Q_0$ 
and denote  $M=Q*_B Q_0$. 

$1^\circ$ If $u\in \Cal N_{Q_0}(B)$ is such that $E_B(u^k)=0$, $\forall 1\leq k\leq n-1$, and we let $\Cal U_0=\{1, u, ..., u^{n-1}\}$, $\varepsilon=2\sqrt{n-1}/n$, 
then $Q\subset M$ is $(\Cal U_0, \varepsilon)$-compressible relative to $B$. 

$2^\circ$ If  there exist $u\in \Cal N_{Q_0}(B)$,  
$v\in \Cal N_Q(B)$,  with $\tau(u)=0$ and $E_B(v^k)=0$, $\forall k\neq 0$, then $Q\subset M$ is compressible relative to $B$.   
If in addition $B \subset M$ is compressible  $($for instance, if $B$ is matrix embeddable$)$, then $Q \subset M$ is compressible.  
\endproclaim
\noindent
{\it Proof}. $1^\circ$ Let $x\in (Q\otimes \Bbb M_K(\Bbb C))_1$ be an element with $0$ expectation onto $B\otimes \Bbb M_K(\Bbb C)$.   

If the condition $1^\circ$ holds, then the set $\{u^kxu^{-k} \mid  k=0, 2, ...,n-1\}\subset (M)_1$ 
is L-free in the sense of (3.1 in [PV14]). By (3.4 in [PV14]), this set can be dilated to a set of L-free unitaries $V_0, V_2, ..., V_{n-1}$ in a larger II$_1$ factor 
$\tilde{M}$ and thus, by the Kesten-type norm estimate 
in [AO77], one has 
$$
\|\frac{1}{n}\sum_{k=0}^{n-1} u^kxu^{-k}\| \leq \|\frac{1}{n}\sum_{k=0}^{n-1} V_k\|\leq 2\sqrt{n-1}/n.
$$ 

For an arbitrary $x\in (Q\otimes \Bbb M_K(\Bbb C))_1$, one writes $x=x'+b$ with $x'$ having $0$ expectation onto $B\otimes \Bbb M_K(\Bbb C)$ 
and $b\in (B\otimes \Bbb M_K(\Bbb C))_1$. The first part shows that 
$\|\frac{1}{n}\sum_{k=0}^{n-1} u^kx'u^{-k}\|\leq 2\sqrt{n-1}/n$, while the fact that $u$ normalizes $B$ insures that 
$\frac{1}{n}\sum_{k=0}^{n-1} u^kbu^{-k}$ stays within $B$. Thus, $\frac{1}{n}\sum_{k=0}^{n-1} u^kxu^{-k}$ is 
$2\sqrt{n-1}/n$-close to $\frac{1}{n}\sum_{k=0}^{n-1} u^kbu^{-k}\in B$. 

$2^\circ$.  The given conditions imply that $u_0=uvu^*$ 
is a Haar unitary normalizing $B$ with the property that the set $\{u_0^kxu_0^{-k}\mid k=1, 2, ...\}$ is L-free, and by the argument in $1^\circ$ above it follows that 
$Q\subset M$ is $(\{u_0^k\mid 0\leq k \leq n-1\}, 2\sqrt{n-1}/n)$-compressible  relative to $B$, $\forall n$. Thus, $Q \subset M$ is compressible  
relative to $B$. 
\hfill $\square$

\proclaim{1.5. Corollary} $1^\circ$ Let $Q, Q_0$ be tracial von Neumann algebras, with $Q$ diffuse and $Q_0$ containing a trace $0$ unitary. Then 
$Q\subset M=Q*Q_0$ is compressible . 
\vskip.05in

$2^\circ$ Let $\Gamma_0, H$ be discrete groups, with $|H|\geq 2$ and $|\Gamma_0|=\infty$. Let $\Gamma_0* H \curvearrowright (B, \tau)$ 
be a trace preserving action. Then the inclusion $B\rtimes \Gamma_0=:Q\subset  M:=B \rtimes (\Gamma_0*H)$ is compressible  relative to $B$. 

\vskip.05in 
$3^\circ$ Let $Q\subset P$ be an  inclusion of $\text{\rm II}_1$ factors 
with finite index, $1< [P:Q]<\infty$,  and $Q\subset P\subset P_1 \subset ... \nearrow P_\infty$ its Jones tower 
and enveloping factor. Let $N_0$ be a tracial von Neumann algebra containing a trace $0$ unitary. 
Then $P_\infty \subset P_\infty *_P (P\overline{\otimes} N_0)$ is compressible  relative to $P$. 

\vskip.05in 
$4^\circ$ Let $M$ be an ultraproduct $\text{\rm II}_1$ factor. Then any separable W$^*$-subalgebra $Q\subset M$ is compressible. 

\vskip.05in 
$5^\circ$ Let $M$ be a $\text{\rm II}_1$ factor of the form $M=N*N_0$ with $N, N_0$ tracial diffuse W$^*$-algebras. Given any 
diffuse W$^*$-subalgebra $B\subset M$, its relative commutant in the ultrapower of $M$ over a free 
ultrafilter $\omega$ on $\Bbb N$, $B'\cap M^\omega$, is compressible. in $M^\omega$. 

\endproclaim 
\noindent
{\it Proof}.  $1^\circ$ If $v\in \Cal U(Q_0)$ has trace $0$, then $vQv^*$ is free independent to $Q$ and diffuse. So $vQv^*$  
contains Haar unitaries $u$. For any such $u$  Proposition 1.4 applies to get that $Q$ is $2\sqrt{n-1}/n$-compressible, $\forall n\geq 2$, 
implying that $Q$ is compressible. 

$2^\circ$ This is immediate from Proposition 1.4. 

$3^\circ$ Let $v\in N_0$ be a trace $0$ unitary and denote $N=P'\cap P_\infty$ (note that this is a II$_1$ W$^*$-algebra, so in particular it 
is diffuse). Note that $N_0=vNv^*$ is free independent to $P_\infty$ relative to $P$, i.e., any alternating word with letters in $P_\infty\ominus P$ 
and $N_0\ominus \Bbb C1$ has trace $0$.  Thus, Proposition 1.4 applies to get that $P_\infty$ is compressible relative to $P$ with 
compressing set of unitaries $\Cal U(N_0)$. 

$4^\circ$ By ([P92], see also [P13]) there exists a diffuse abelian W$^*$-subalgebra $A\subset M$ that's free independent to $Q$, 
so Proposition 1.4 applies to get that $Q\subset M$ is compressible. 

$5^\circ$ Let $A\subset B$ be an abelian diffuse W$^*$-subalgebra. Since $B'\cap M^\omega \subset A'\cap M^\omega$, 
it is sufficient to prove that $A'\cap M^\omega$ is compressible. 
Since any two embeddings in $M^\omega$ of the separable abelian diffuse W$^*$-algebra are unitarily conjugate, 
one may assume $A\subset N$. By [HI23], $A'\cap M^\omega$ is free independent to $A_0'\cap M^\omega$, for any diffuse abelian W$^*$-algebra $A_0\subset N_0$. 
So in particular, there exist Haar unitaries $u\in A_0'\cap M^\omega$ so that $\{u^n\}_n$ is free independent to $A'\cap M^\omega$. 
Proposition 1.4 then applies to get that $A'\cap M^\omega$ is compressible in $M^\omega$. 

\hfill $\square$

\proclaim{1.6. Proposition} Assume $M=N *A$, with $N\neq \Bbb C1$ 
and $A$ abelian diffuse. Then there exists an increasing sequence of compressible  W$^*$-subalgebras $M_n \subset M$ 
such that $\overline{\cup_n M_n}^w=M$. If moreover $N$ is diffuse, then $M_n \subset M$ can be taken irreducible with spectral gap, $\forall n$. 
\endproclaim
\noindent
{\it Proof}. By (Corollary 7.4 in [P90]), $M$ is a non-Gamma II$_1$ factor. Let $A_n\subset A$ 
be an increasing sequence of dyadic partitions exhausting $A$ and define $M_n=N*A_n$. 
We clearly have $M_n \nearrow M$. Also, for each $n$ we have $M=M_n*_{A_n} A$. 
By Proposition 1.4.2$^\circ$, this implies that $M_n\subset M$ is compressible  relative to $A_n$. 
Since $A_n \subset M$ is also compressible, by 1.3.4$^\circ$ it follows that $M_n \subset M$ is compressible. 

If we now assume $N$ is diffuse, then by [P90] again, the algebras $M_n$ 
are non-Gamma II$_1$ factors with $M_n\subset M$ irreducible with spectral gap, $\forall n$.  

\hfill $\square$

\vskip.05in 

We denote as usual by $L\Bbb F_t, 1 < t \leq \infty$, the {\it interpolated free group factors} as introduced in ([Dy94], [R94]). 

\proclaim{1.7. Corollary} If $M=L\Bbb F_t, 1< t \leq \infty$, then $M$ admits an increasing sequence of compressible  subfactors 
$M_n \subset M$ 
such that $\overline{\cup_n M_n}^w=M$. If in addition $t\geq 2$, then $M_n\subset M$ can be taken irreducible with spectral gap. 
\endproclaim
\noindent
{\it Proof}. By [Dy93], one can write each interpolated free group factor $M=L\Bbb F_t$ as $N*A$ for some $N\neq \Bbb C1$ and $A$ abelian diffuse.  
Moreover, if $t\geq 2$ then $N$ can be taken diffuse.  The statement then follows from Proposition 1.6. 

\hfill $\square$

\vskip.05in 
\noindent
{\bf 1.8. Remarks}. $1^\circ$ We do not have any example of a tracial W$^*$-inclusion $Q\subset M$ with $Q$ diffuse that's compressible 
but for which there exists no diffuse abelian subalgebra $A\subset M$ such that $A$ is free independent to $Q$.  

$2^\circ$ It would be interesting to know whether there exist or not tracial W$^*$-inclusions $Q\subset M$ 
that are not compressible but for which the following {\it weak compressibility} condition holds true:  

\vskip.05in

{\it $\forall \varepsilon >0$, $\exists u_1, ..., u_n \in \Cal U(M)$ such that $\|\frac{1}{n} \sum_{j=1}^n u_j xu_j^*-\tau(x)1\| \leq \varepsilon$, 
$\forall x\in (Q)_1$}.

\vskip.05in

Also, does the existence of a diffuse $Q\subset M$ satisfying this  weaker condition imply $M$ is non-amenable? 

$3^\circ$ It would be interesting to have a ``lucrative'' sufficient condition for an inclusion of infinite groups $H\subset G$ 
that insures $LH \subset LG$ is compressible, respectively weakly compressible 
(and that would of course not involve existence of an infinite $H_0\subset G$ that's free  to $H$!). 


\heading 2. Combining compressibility with weak-coarseness    
\endheading

Recall that the coarse  Hilbert $\Cal N-\Cal M$ bimodule over W$^*$-algebras $\Cal N, \Cal M$ 
is defined as the tensor product of their standard representations, $L^2\Cal N \overline{\otimes} L^2\Cal M$, with $\Cal N$ 
acting on $L^2\Cal N$ from the left and $\Cal M$ on $L^2\Cal M$ from the right. 
More generally, we will refer to any sub-bimodule of $(L^2\Cal N \overline{\otimes} L^2\Cal M)^{\oplus \infty}$ as a coarse 
$\Cal N-\Cal M$ bimodule (see e.g. [P86], [P18]). 

It is trivial to see that any $\Cal N-\Cal M$ 
Hilbert bimodule $_\Cal N\Cal H_\Cal M$ contains a maximal coarse sub-bimodule, i.e., a sub-bimodule $\Cal H' \subset \Cal H$ 
such that $\Cal H'\subset (L^2\Cal N\overline{\otimes} L^2\Cal M)^{\oplus\infty}$ and $\Cal H\ominus \Cal H'$ contains 
no non-zero coarse sub-bimodule. We denote this sub-bimodule of $_\Cal N\Cal H_\Cal M$ by $(_\Cal N\Cal H_\Cal M)^{co}$ and call it the {\it coarse part} of $\Cal H$. 

Recall also from (2.1.5 in [P86]; see also [A95], [AP17]) 
that a Hilbert $\Cal N -\Cal M$ bimodule $\Cal H$ is weakly contained in a $\Cal N -\Cal M$ bimodule $\Cal H'$, written $\Cal H\prec \Cal H'$,  
if $\Cal H$ is in the closure of $\Cal H'$ in the Fell-type topology on bimodules (see 2.1.1 in [P86]). This simply means that 
$\Cal H$ can be ``simulated locally'' inside $\Cal H'$. This condition is easily seen to be equivalent to the norm on $\Cal N \vee_{Alg} \Cal M^{op}$ in 
its representation on $\Cal H$ 
being majorized by the norm in its representation on $\Cal H'$. 

\vskip .05in

\noindent
{\bf 2.1. Definition.} A Hilbert $\Cal N-\Cal M$ bimodule $_\Cal N\Cal H_\Cal M$ over the W$^*$-algebras $\Cal N, \Cal M$ is  {\it weakly coarse} if  
$_\Cal N\Cal H_\Cal M$ is weakly contained in the coarse $\Cal N-\Cal M$ bimodule. 

By the remarks above, this amounts to $\|x\|_{\Cal B(\Cal H)}\leq \|x\|_{\Cal N\overline{\otimes} \Cal M}$, for all $x$ in $\Cal N \otimes \Cal M^{op}$. 

\vskip.05in

Weak coarseness is preserved by taking sub-bimodules and direct sums. A coarse bimodule is of course weakly coarse. 
A well known result due to Effros and Lance ([EL77]) shows that 
any Hilbert bimodule in which one of the algebras involved is AFD is weakly coarse. But in general, such a bimodule is not coarse. 
For instance, if $\Cal N \subset \Cal M$ is a proper W$^*$-inclusion with $\Cal N$ AFD and if either $\Cal N'\cap \Cal M=\Bbb C$, or 
if there exists a normal expectation of $\Cal M$ onto $\Cal N$, then $_{\Cal N}L^2\Cal M_{\Cal M}$ has coarse part equal to $0$.  

We record this example below, together with several other interesting examples of weakly coarse bimodules 
(for the definition of tensor product, or ``Connes fusion'', of Hilbert bimodules, mentioned in part $3^\circ$ see e.g., $1.3.1$ in [P86]).

\proclaim{2.2. Proposition} $1^\circ$ $(\text{\rm [EL77]})$ Any Hilbert $\Cal N-\Cal M$ bimodule where $\Cal N$ or $\Cal M$ is AFD, is weakly coarse. 

$2^\circ$ $\text{\rm ([IT23], [DiP23])}$ If $M$ is an interpolated free group factor, or more generally an amplification of a group factor 
$L\Gamma$ with $\Gamma$ bi-exact, then $_M(L^2(M^\omega) \ominus L^2M)_M$ is  weakly coarse. 

$3^\circ$ Given any bimodules $\Cal H =$ $_N\Cal H_B$, $\Cal K =$ $_B\Cal K_M$, where $B, N, M$ are tracial W$^*$-algebras and $B$ is AFD, 
the tensor product $_N(\Cal H\overline{\otimes}_{\Cal B} \Cal K)_M$ is weakly coarse. 

\endproclaim
\noindent
{\it Proof} Part 1$^\circ$ is in [EL77]. Part $2^\circ$ is shown in (Lemma 4.4 of [IT23]) in the case $M=L\Bbb F_n^t$, $t>0$, and in ([DiP23]) 
for $M=L\Gamma^t$, $t>0$, when $\Gamma$ is an arbitrary bi-exact group. 

Part 3$^\circ$ is a folklore result and follows easily by noticing that if $B_0$ is finite dimensional (or merely atomic) 
then $_N(\Cal H\overline{\otimes}_{B_0} \Cal K)_M$ is coarse and that  if $B_n \nearrow B$ are finite dimensional 
then the $N-M$ bimodules $\Cal H'_n:=\Cal H\overline{\otimes}_{B_n} \Cal K$ converge to $\Cal H\overline{\otimes}_{B} \Cal K$. 
\hfill $\square$

\vskip.05in

We'll now prove that  compressibility 
forces a weakly coarse bimodule to always contain a copy of the coarse bimodule, and in certain situations to be  ``totally'' coarse.

\proclaim{2.3. Theorem} Let $_\Cal N \Cal H_\Cal M$ be a Hilbert bimodule over the W$^*$-algebras $\Cal N, \Cal M$. Let $Q, B \subset \Cal M$ 
be W$^*$-subalgebras. Assume $_{\Cal N}\Cal H_{\Cal M}$ is weakly coarse and 
$Q\subset \Cal M$ is compressible relative to $B$.

Then $\Cal N\vee B^{op}  = ({\Cal M^{op}}' \cap (\Cal N \vee Q^{op}))\vee B^{op}$. 
\endproclaim
\noindent
{\it Proof}. Since $\Cal N \vee B^{op}\subset ({\Cal M^{op}}' \cap (\Cal N \vee Q^{op})) \vee B^{op}$, all we need to prove is that ${\Cal M^{op}}' \cap (\Cal N \vee Q^{op}) \subset \Cal N\vee B^{op}$, 
in other words that ${\Cal M^{op}}' \cap (\Cal N \vee Q^{op})$ is contained in the $so$-closure of the $^*$-algebra $\Cal N \vee_{alg} B^{op}$. 

To this end, it is sufficient to show that given any $y \in ({\Cal M^{op}}' \cap (\Cal N \vee Q^{op}))_1$, any finite 
 $F\subset (\Cal H)_1$ and any $\varepsilon >0$,   
there exists $X^0 \in (\Cal N \vee_{alg} B^{op})_1$ such that $\|(y-X^0)(\xi)\|_{\Cal H}\leq \varepsilon$, 
$\forall \xi\in F$. 

By the  compressibility of $Q\subset \Cal M$ relative to $B$ and the weak coarseness of $_{\Cal N}\Cal H_{\Cal M}$ (used here in its equivalent 
formulation stated in the last part of 2.1), it follows that 
there exist unitary elements $u_1, ..., u_n\in \Cal M^{op}$ such that  
for any $X\in (\Cal N \vee_{alg} Q^{op})_1$ there exists $X^0\in (\Cal N \vee_{alg} B^{op})_1$ with the property that 
$$
\| \frac{1}{n}\sum_{j=1}^n u_j X u_j^* - X^0 \| \leq \varepsilon/2. \tag 2.3.1  
$$

Since  $y$ lies in the W$^*$-algebra generated by the $^*$-algebra $\Cal N \vee_{alg} Q^{op}$, 
by using Kaplansky's density theorem it follows that there exists $X \in (\Cal N \vee_{alg} Q^{op})_1$ such that 
for all $\xi\in F$ and all $1\leq j \leq n$ we have 
$$
\|(y-X)(u_j^*(\xi))\|_{\Cal H} \leq \varepsilon/2 . \tag 2.3.2 
$$ 
By applying the unitary $u_j$ to the vector $(y-X)(u_j^*(\xi))\in \Cal H$ in $(2.3.2)$ 
and taking into account that $u_j\in \Cal M^{op}$ commutes with $y \in {\Cal M^{op}}'$, 
we thus get for each $\xi\in F$ the estimate
$$
\|(y-u_jXu_j^*)(\xi)\|_{\Cal H}=\|u_j(y-X)(u_j^*(\xi))\|_{\Cal H}  \tag 2.3.3
$$
$$
= \|(y-X)(u_j^*(\xi))\|_{\Cal H}\leq \varepsilon/2,  
$$
which by summing up over $j$ and dividing by $n$ gives
$$
\|(y-\frac{1}{n}\sum_{j=1}^nu_jXu_j^*)(\xi)\|_{\Cal H} \leq \frac{1}{n}\sum_{j=1}^n \|(y-u_jXu_j^*)(\xi)\|_{\Cal H} \leq \varepsilon/2 \tag 2.3.4 
$$

But by $(2.3.1)$, there exists $X^0\in (\Cal N \vee_{alg} B^{op})_1$ such that  
$$
\| \frac{1}{n}\sum_{j=1}^n u_j X u_j^* - X^0 \| \leq \varepsilon/2. \tag 2.3.5 
$$
Thus, by combining $(2.3.4)$ and $(2.3.5)$ we finally get 
$$
\|(y-X^0)(\xi)\|_{\Cal H} \leq \varepsilon/2 + \varepsilon/2=\varepsilon, 
$$
for all $\xi \in F$. 
 
\hfill $\square$

\proclaim{2.4. Theorem} Let $_\Cal N \Cal H_\Cal M$ be a weakly coarse Hilbert bimodule over W$^*$-algebras $\Cal N, \Cal M$. Let $Q, B \subset \Cal M$ 
be W$^*$-subalgebras. 

$1^\circ$  If $\Cal U_0 \subset \Cal U(\Cal M)$ finite, $\varepsilon >0$ are such that 
$Q \subset \Cal M$ is $(\Cal U_0, \varepsilon)$-compressible relative to $B$,  
then $\Cal N \vee Q^{op}$ is $(\Cal U_0^{op}, \varepsilon)$-compressible relative to $\Cal N \vee B^{op}$. 

$2^\circ$ If $Q\subset \Cal M$ is compressible, then $\Cal N \vee Q^{op}$ is compressible relative to $\Cal N\vee B^{op}$. 
\endproclaim
\noindent
{\it Proof}. $1^\circ$ Let $\Cal U_0^{op}=\{u_1, ..., u_n\}$. 
Let $y\in (\Cal N\vee Q^{op})_1$.  
Since  $y$ lies in the von Neumann  
algebra generated by the $^*$-algebra $\Cal N \vee_{alg} Q^{op}$, 
by using Kaplansky's density theorem it follows that for any $F\subset (\Cal H)_1$ 
there exists $X_F \in (\Cal N \vee_{alg} Q^{op})_1$ such that 
for all $\xi\in F$ and all $1\leq j \leq n$ we have 
$$
\|(y-X_F)(u_j^*(\xi))\|_{\Cal H} \leq 1/|F|. \tag 2.4.1 
$$ 
By applying the unitary $u_j$ to the vector $(y-X_F)(u_j^*(\xi))\in \Cal H$ in $(2.4.1)$  
we thus get for each $\xi\in F$ and $1\leq j \leq n$, the estimate
$$
\|(u_jyu_j^*-u_jX_Fu_j^*)(\xi)\|_{\Cal H}=\|u_j(y-X_F)(u_j^*(\xi))\|_{\Cal H} \tag 2.4.2 
$$
$$
= \|(y-X_F)(u_j^*(\xi))\|_{\Cal H}\leq 1/|F|, 
$$
which by summing up over $j$ and dividing by $n$ gives
$$
\|(\frac{1}{n}\sum_{j=1}^nu_jyu_j^*-\frac{1}{n}\sum_{j=1}^nu_jX_Fu_j^*)(\xi)\|_{\Cal H}  \tag 2.4.3 
$$
$$
\leq \frac{1}{n}\sum_{j=1}^n \|(u_jyu_j^*-u_jXu_j^*)(\xi)\|_{\Cal H} \leq 1/|F|. 
$$

On the other hand, by the $(\Cal U_0, \varepsilon)$-compressibility of $Q\subset \Cal M$ relative to $B$,  
for each $X_F\in (\Cal N \vee_{alg} Q^{op})_1$ there exists $X^0_F\in (\Cal N \vee_{alg} B^{op})_1$ with the property that 
$$
\| \frac{1}{n}\sum_{j=1}^n u_j X_F u_j^* - X_F^0 \| \leq \varepsilon. \tag 2.4.4  
$$
Thus, by combining $(2.4.3)$ and $(2.4.4)$ we finally get 
$$ 
\|(\frac{1}{n}\sum_{j=1}^nu_jyu_j^*-X_F^0)(\xi)\|_{\Cal H} \leq \varepsilon +  1/|F| \tag 2.4.5
$$
for all $\xi \in F$. Taking a weak limit point $b$ of the net $\{X^0_F\}_F \subset (\Cal N \vee_{alg} B^{op})_1$, 
indexed over the directed set all finite subsets $F\subset (\Cal H)_1$, with the inclusion order, 
by the compactness of  $(\Cal N \vee B^{op})_1$ and the upper semicontinuity of the operator norm with respect to the weak topology, it follows   
that $b\in (\Cal N \vee B^{op})_1$ and that by $(2.4.5)$ it satisfies $\|\frac{1}{n}\sum_{j=1}^nu_jyu_j^*-b\|\leq \varepsilon$. 

Part $2^\circ$ is an immediate consequence of $1^\circ$. 
\hfill $\square$ 

\vskip.1in

\proclaim{2.5. Theorem} Let $_\Cal N \Cal H_M$ be a weakly coarse Hilbert bimodule over W$^*$-algebras $\Cal N, M$. 
Assume  $B\subset Q \subset M$ 
are W$^*$-subalgebras such that $Q \subset M$ is compressible  relative to $B$ with compressing set of unitaries $\Cal U\subset \Cal U(M)$ 
satisfying the following conditions: $(a)$ $[\Cal U, B]=0$; 
$(b)$ there exists a normal conditional expectation $E_B: Q\rightarrow B$ such that $E_B(uxu^*)=E_B(x)$, $\forall x\in Q$, $u\in \Cal U$. 
Then we have:

\vskip.05in

$1^\circ$ There exists a unique normal conditional expectation $\Phi: \Cal N \vee Q^{op}\rightarrow \Cal N \vee B^{op}$ 
such that $\Phi(x y^{op})=xE_B(y)^{op}$ for all $x\in \Cal N, y\in Q$. 

$2^\circ$ Denote by $\Cal N_0$ the von Neumann algebra generated by the normalizer of $\Cal N$ in $(Q^{op})'\cap \Cal B(\Cal H)$ 
and $Q_0:=\Cal N_Q(B)''\subset Q$. Then the support projection 
$p_{\phi}$ of $\Phi$ $($i.e., the largest projection $p\in \Cal P(\Cal N\vee Q^{op})$ with the property that $x\in (p\Cal N \vee Q^{op}p)_+$ 
$\Phi(x)=0$ implies $x=0)$ satisfies $p_{\Phi}\in (\Cal N_0 \vee Q_0^{op})'\cap (\Cal N \vee Q^{op})$. In particular, if $\Cal N_Q(B)''=Q$, 
then $p_\Phi \in \Cal N_0'\cap \Cal Z(\Cal N \vee Q^{op})$. 

$3^\circ$ Asume $B=\Bbb C$. For any non-zero $\xi\in \Cal H$, the $\Cal N-Q$ Hilbert bimodule $\Cal H_\xi=\overline{sp} \Cal N \xi \Cal U Q$ contains 
a coarse $\Cal N-Q$ Hilbert sub-bimodule: 
there exists $0\neq \Cal H_0\subset \Cal H_\xi$ such that $\Cal N \Cal H_0 Q\subset \Cal H_0$ and $_\Cal N{\Cal H_0}_Q \simeq L^2(\Cal N z)\overline{\otimes} L^2Q$, 
where $z\in \Cal P(\Cal Z(\Cal N))$ is the support projection of the representation $\Cal N\ni x\mapsto xp_{\Cal H_0}\in \Cal B(\Cal H_0)$. 
\endproclaim
\noindent
{\it Proof}. $1^\circ$ Let $E:Q^{op} \rightarrow B^{op}$ be defined by $E(x^{op})=E_B(x)^{op}$, $\forall x\in Q$.  
Denote by  $\Phi_0: \Cal N \otimes_{\min} Q^{op} \rightarrow \Cal N\otimes_{\min} B^{op}$ 
the conditional expectation $id_{\Cal N}\otimes E: \Cal N\otimes_{\min} Q^{op} \rightarrow \Cal N\otimes_{\min} B^{op}$. 

Note that by the compressibility of $Q\subset \Cal M$ relative to $B$ and the assumptions 
on $E$, we have the following: 

\vskip.05in 

{\it Fact}. For any $\varepsilon_0>0$ there 
exists $\Cal U_0\subset \Cal U^{op}$ finite such that 
$$
\|\frac{1}{|\Cal U_0|} \sum_{u\in \Cal U_0} uxu^* - \Phi_0(x)\| \leq \varepsilon_0, \forall x\in (\Cal N \otimes_{\min} Q^{op})_1 \tag 2.5.1
$$ 
Indeed, by the compressibility assumption there exists $\Cal U_0\subset \Cal U^{op}$ finite such that for each 
$x\in (\Cal N \otimes_{\min} Q^{op})_1$ there exists $b(x)\in (\Cal N \otimes_{\min} B^{op})_1$ satisfying   
$$
\|\frac{1}{|\Cal U_0|} \sum_{u\in \Cal U_0} uxu^*-b(x)\| \leq \varepsilon_0/2.
$$ 
By applying $\Phi_0=id_{\Cal N}\otimes E$, which is 
a conditional expectation and thus contractible and using that $E(uyu^*)=E(y)$, $\forall y\in Q^{op}$, $u\in \Cal U_0$, this implies     
$\|\Phi_0(x)-b(x)\|\leq \varepsilon_0/2$ and $(2.5.1)$ follows by the triangle inequality.  

\vskip.05in

Now, if we identify  $\Cal N\otimes_{\min} Q^{op}$ with $\Cal N\vee_{C^*} Q^{op}$ and $\Cal N\otimes_{\min} B^{op}$ 
with $\Cal N\vee_{C^*} B^{op}$,  
then in order to prove the statement 
we need to show that $\Phi_0$ extends to an so-continuous map $\Phi$ from $\Cal N\vee Q^{op}$ to $\Cal N \vee B^{op}$.

By Kaplansky's density theorem, it is sufficient to show that if $T_i=\sum_k x_k^i y^i_k \in (\Cal N\vee_{alg} Q^{op})_1$ 
is a net converging in the so-topology to some $T\in \Cal B(\Cal H)$, 
then $\Phi_0(T_i)$ is so-convergent. This amounts to showing that given any $F_0\subset (\Cal H)_1$ finite and $\varepsilon>0$, there exists a finite set of indices 
$J$ such that if $i, j \not\in J$ then $\|(\Phi_0(T_i)-\Phi_0(T_j))(\xi)\|<\varepsilon$, $\forall \xi\in F_0$. 

By the above {\it Fact}, there  exists $\Cal U_0\subset \Cal U^{op}$ finite  such that for any $X=\sum_k x_k \otimes y_k \in (\Cal N \otimes Q^{op})_1 \subset (\Cal N \otimes_{min} \Cal M^{op})_1$ we have 
$$
\|\frac{1}{|\Cal U_0|} \sum_{u\in \Cal U_0} uXu^* - \Phi_0 (X)\|< \varepsilon/3. \tag 2.5.2
$$
By the so-convergence of the net $\{T_i\}_{i\in I}$, there exists 
a finite set $J\subset I$ such that if we denote $F=\{u^*(\xi)\mid u\in \Cal U_0, \xi \in F_0\}$, then $\|(T_i-T_j)(\eta)\|< \varepsilon/3$, for all $i, j \in I \setminus J$ and all $\eta \in F$. 
Thus 
$$
\|u(T_i-T_j)u^*(\xi)\| = \|(T_i-T_j)(u^*(\xi))\| < \varepsilon/3, \forall u\in \Cal U_0, \xi\in F_0 \tag 2.5.3 
$$
Since 
each $u\in \Cal U_0$ commutes with $\Phi_0(T_i)$,  $\forall i\in I$ (because $[\Cal U, B]=0$), it follows that if  $\xi\in F_0$ and $i, j \in I\setminus J$, then 
by using $(2.5.2)$ and $(2.5.3)$ we get   
$$
\|(\Phi_0(T_i)-\Phi_0(T_j))(\xi)\|  \tag 2.5.4 
$$
$$
\leq \|\Phi_0(T_i)-\frac{1}{|\Cal U_0|} \sum_u uT_iu^*\| 
+ \frac{1}{|\Cal U_0|} \sum_u \|u(T_i-T_j)u^*(\xi)\| 
$$
$$
+ \|\frac{1}{|\Cal U_0|} \sum_u uT_ju^* - \Phi_0(T_j)\| 
\leq \varepsilon/3 + \varepsilon/3 + \varepsilon/3=\varepsilon. 
$$

$2^\circ$ To show that $p=p_\Phi$ is in $(\Cal N_0 \cup Q^{op}_0)'\cap (\Cal N \vee Q^{op})$ 
notice first that the $\Cal N$-bimodularity of $\Phi$ implies 
$p\in \Cal N'\cap (\Cal N\vee Q^{op})$. 

If now we take $u\in \Cal N_Q(B)^{op}$, then by the formula of $\Phi_0$ we have $\Phi(uXu^*)=\Phi_0(u X u^*)=u\Phi_0(X)u^*=u\Phi(X)u^*$ for all 
$X\in \Cal N \vee_{Alg} Q^{op}$, which by the so-continuity of $\Phi$ implies $\Phi(uXu^*)=u\Phi(X)u^*$, $\forall X\in \Cal N \vee Q^{op}$. This implies that if 
$x\in (p\Cal N \vee Q^{op}p)_+$ satisfies $\Phi(x)=0$ then $\Phi(uxu^*)=0$. Hence,  
$upu^*=p$, $\forall u\in \Cal N_Q(B)^{op}$,  and thus $p\in (\Cal N_Q(B)^{op})' = (Q_0^{op})'$, 
implying that $p\in (\Cal N \cup Q_0^{op})'\cap \Cal N \vee Q^{op}$. 

Similarly, if we take $v$ in the normalizer of $\Cal N$ in $(Q^{op})'\cap \Cal B(\Cal H)$, then 
by the formula of $\Phi_0$ we have $\Phi(vXv^*)=\Phi_0(v X v^*)=v\Phi_0(X)v^*=v\Phi(X)v^*$,  
$\forall X\in \Cal N \vee_{Alg} Q^{op}$, and thus $\forall X\in \Cal N\vee Q^{op}$. By the above argument, this implies $p\in \Cal N_0'$.

This also shows that if $Q_0=Q$, then  
$p\in (\Cal N_0 \cup Q^{op})'\cap \Cal N \vee Q^{op}=\Cal N_0'\cap \Cal Z(\Cal N \vee Q^{op})$.

$3^\circ$  Since $B=\Bbb C$, the map $\Phi: \Cal N \vee Q^{op} \rightarrow \Cal N$ defined in Part $1^\circ$ 
becomes a normal  conditional expectation onto $\Cal N$ and the statement follows immediately from part $2^\circ$ above. 
\hfill $\square$

\vskip.05in

For the next result, we'll say that a subgroup $H$ of a (discrete) group $G$ is compressible if the corresponding inclusion of group 
W$^*$-algebras $LH\subset LG$ is compressible. Also, if $\mycal R$ is a countable measurable measure preserving equivalence relation 
on the standard probability measure space $(X, \mu)$ and $\Cal S\subset \mycal R$ is a subequivalence relation, then we say that $\Cal S\subset \mycal R$ 
is compressible if the associated W$^*$-inclusion $L(\Cal S)\subset L(\mycal R)$ is compressible relative to $L^\infty X$. 

\proclaim{2.6. Corollary} $1^\circ$ If a tracial W$^*$-algebra $M$ contains  
a diffuse compressible subalgebra, then $M$ does not have any AFD direct summand.   

\vskip.05in
$2^\circ$ If a group $G$ has an infinite compressible subgroup $H\subset G$, then $G$ is not amenable.  

\vskip.05in
$3^\circ$ If $\mycal R$ is an ergodic measure preserving equivalence relation on a probability measure space $(X, \mu)$  that contains 
a compressible subequivalence relation $\Cal S$ of type $\text{\rm II}_1$, then $\mycal  R$ is not amenable.   
\endproclaim 
\noindent
{\it Proof}.  $1^\circ$ Assume by contradiction that $M$ has a non-zero  AFD direct summand. By ``cutting'' if necessary  with a central projection, 
we may assume $M$ itself is AFD and that it has a diffuse abelian von Neumann subalgebra $A\subset M$ that's 
compressible. The compressibility of $A\subset M$ with $A$ diffuse implies that $M$  cannot have type I direct summands. 

So $M$ is necessarily of type II$_1$, with $_M L^2M_A$ weakly coarse and  $A^{op}\subset M^{op}$ compressible. 
By  Theorem 2.5, the map $\Cal E: M \vee_{alg} A^{op}\rightarrow M$ defined by $\Cal E(x a^{op})=\tau(a)x$, $x\in M, a\in A$, extends to 
a normal conditional expectation of $\Cal M=M \vee A^{op}$ onto $M\subset \Cal M$. But $\Cal M$ is properly infinite, $M$ is a finite von Neumann subalgebra 
and $M'\cap \Cal M$ has no finite projections of $\Cal M$, implying that there are no normal conditional expectations of $\Cal M$ onto $M$.

$2^\circ$ Since the inclusion $N=L(H)\subset L(\Gamma)=M$ is compressible, part 1$^\circ$ applies to get that $L(\Gamma)$ is not amenable, thus $\Gamma$ 
is not amenable. 

$3^\circ$ If we denote $N=L(\Cal S)\subset L(\Cal R)=M$ and $A=L^\infty(X)\subset N$, then the 
hypothesis states that $N\subset M$ is compressible relative to $A$, with $M$ a II$_1$ factor and $N$ a II$_1$ von 
Neumann algebra that contains the Cartan subalgebra $A$ of $M$. 

If we assume $M\simeq R$, then we claim that  
there exists a II$_1$ subfactor $M_0\subset M$ of index $2$ in $M$ 
such that $N'\cap \langle M, e_{M_0}\rangle =N'\cap M$. To see this, let $\theta\in \text{\rm Aut}(M)$ be a period 2 automorphism 
with the property that $\theta(A)\perp A$. If $\tilde{M}$ denotes the crossed product II$_1$ factor $M\rtimes_\theta \Bbb Z/2\Bbb Z$ 
with $u\in \tilde{M}$ the canonical unitary implementing $\theta$ on $M$, then $uAu^*\perp A$. By (Lemma 2.5 in [P81b]), it follows that the space 
$uM$ is perpendicular to the von Neumann algebra $M_1$ generated by the normalizer of $A$ in $\tilde{M}$. 
Taking into account that $L^2(\tilde{M})=L^2M \oplus L^2(uM)$, we thus get $M_1\subset (uM)^\perp=(\tilde{M}\ominus M)^\perp=M$.  
In particular $A'\cap \tilde{M}\subset M_0\subset M$ and hence $N'\cap \tilde{M}\subset A'\cap \tilde{M}=A\subset N$, 
implying that $N'\cap \tilde{M}\subset N$. 

Thus, if we take $M_0\subset M$ to be the fixed point algebra $\{x\in M \mid \theta(x)=x\}$, then $(M\subset \tilde{M})=(M\subset \langle M, e_{M_0}\rangle)$ 
and we indeed have  $N'\cap \langle M, e_{M_0}\rangle =N'\cap M$.  
But this implies $M_0\vee N^{op}=M\vee N^{op}$. Thus, ${M^{op}}'\cap (M_0\vee N^{op})=M$. 

On the other hand, since the bimodule $_{M_0}L^2M_N$ is weakly coarse, we can apply Theorem 2.3 to get ${M^{op}}' \cap (M_0\vee N^{op}) = M_0$, a contradiction.

\hfill $\square$


\vskip.05in

We now combine the compressibility condition with a weak coarseness assumption on correlation bimodules 
involving subalgebras of a tracial W$^*$-algebra $M$. Thus, using the above theorems, we show that    
if $Q\subset M$ is compressible with compressible set of unitaries lying in some $P\subset M$ containing $Q$,  
then any weakly coarse bimodule $_NL^2M_P$ has a  non-zero part that's $N-Q$ coarse.  

It is easy to see that if $N, Q$ are W$^*$-subalgebras of the tracial W$^*$-algebra $M$, then the coarse part of the correlation bimodule $_NL^2M_Q$ 
is $\Cal N_M(N)-\Cal N_M(Q)$ invariant. (N.B.: this has already been mentioned in the proof of 2.5.2$^\circ$.)

In fact, the same invariance holds  true with respect to the quasi-normalizers of $N, Q$. 
Recall from (1.4.2 in [P01]) that if $B\subset M$ is a tracial W$^*$-inclusion, then the quasi-normalizer of $B$ in $M$, denoted $q\Cal N_M(B)$, is the set 
of all $x\in M$ with the property that  $\text{\rm sp}(BxB)$ is finitely generated  both as a left $B$-module and as a right $B$-module.   
It is easily seen to be a $^*$-subalgebra of $M$.  

\proclaim{2.7. Lemma} Let $M$ be a tracial W$^*$-algebra and $N, Q\subset M$ be W$^*$-subalgebras. Denoting  
$\Cal H=(_NL^2M_Q)^{co}$,  
$\tilde{N}=\overline{q\Cal N_M(N)}^w$,  
$\tilde{Q}=\overline{q\Cal N_M(Q)}^w$, we have $\tilde{N} \Cal H \tilde{Q}=\Cal H$. 
\endproclaim 
\noindent
{\it Proof}.  Note first that $\Cal H=(_NL^2M_Q)^{co}$ coincides with the set of vectors $\xi \in L^2M$ with the property that  
the positive functional $\varphi_\xi$ 
on $N\vee_{Alg} Q^{op}\subset \Cal B(L^2M)$ defined by  
$\varphi_\xi(xy^{op})=\langle x\xi y, \xi\rangle=\tau(\xi^*x\xi y)$, $\forall x\in N, y\in Q$, is majorized on the restriction on the image of this algebra 
via the quotient map into $N\otimes_{min}Q^{op}$
by a positive element in $L^1(N\overline{\otimes} Q^{op}, \tau_{|N}\otimes \tau_{|Q^{op}})$. 

Also,  if we denote by 
$\Cal H_0=(_NL^2M_Q)^{co}_0$ the set of vectors $\xi\in L^2M$ with the property that $\varphi_\xi$  is majorized by 
$C_{\xi} \tau_{|N}\otimes \tau_{|Q^{op}}$, for some $C_{\xi}>0$, then $\Cal H_0$ is a dense vector subspace of $\Cal H$.

Both $\Cal H_0, \Cal H$ are obviously $\Cal N_M(N)-\Cal N_M(Q)$ invariant. In particular, 
$\Cal H_0, \Cal H$ are $N-Q$ invariant and $(N'\cap M)-(Q'\cap M)$ invariant. 

This description of $\Cal H_0$ easily implies that if $N_0\subset N$ is a 
W$^*$-subalgebra with finite Pimsner-Popa index ([PiP83]), then  $(_{N_0}L^2M_Q)^{co}_0=(_NL^2M_Q)^{co}_0$. 
This implies that if $\psi : p_0Np_0 \rightarrow pNp$ is a unital 1-to-1 $^*$-morphism with image having finite index and $v\in M$ is a 
``quasi-normalizing'' partial isomorphism satisfying $xv = v\psi(x)$, $\forall x\in p_0Np_0$, then $v\Cal H_0 \subset \Cal H_0$ 
(see e.g. the proof  of Theorem 2.1 in [P03]). 
Then note that any $x\in q\Cal N_M(N)$ can be written as a finite sum of the form 
$\sum_i x_i x_i' v_i y_iy_i'$, where $v_i$ are such quasi-normalizing partial isometries and $x_i, y_i\in N, x_i', y_i'\in N'\cap M$. 
Since each term of this summation leaves $\Cal H_0$ invariant, it follows that $x\Cal H_0\subset \Cal H_0$, thus $x\Cal H\subset \Cal H$ 
as well. 

The fact that $\Cal H_0, \Cal H$ are right $Q$-invariant as well is similar.  
 We leave the details as an exercise. 
\hfill $\square$

\vskip.05in 
Recall from ([P01]) that a W$^*$-subalgebra $B\subset M$ is {\it quasi-regular} in $M$ if $q\Cal N_M(B)$ is weakly dense in $M$.

\proclaim{2.8. Corollary} Let $M$ be a tracial W$^*$-algebra and  $P \subset M$ a diffuse quasi-regular W$^*$-subalgebra. 
Then $P$ cannot be compressible. 
\endproclaim
\noindent
{\it Proof}. Assume $P$ is compressible and let $A\subset P$ be a MASA. Note that $P$ diffuse implies $A$ diffuse. Since $A$ is AFD, $_AL^2M_P$ is weakly coarse, 
so  Theorem 2.5 implies that the coarse part $\Cal H=(_AL^2M_P)^{co}$ is non-zero. 
Since $P$ is quasi-regular,  by Lemma 2.7 it follows that  $\Cal H=L^2M$. But this is a contradiction, since $_AL^2M_P$ contains $_AL^2P_P$. 
\hfill $\square$

\vskip.05in 

Following [P18], if  $_NL^2M_Q$ is coarse for some subalgebras $N, Q\subset M$, 
i.e., $_NL^2M_Q$ $\subset (L^2N\overline{\otimes} L^2Q)^{\oplus\infty}$, then we also say that $N$ {\it is coarse to} $Q$.

\proclaim{2.9. Corollary} Let $Q \subset M$ be a compressible tracial W$^*$-inclusion. 
Given any AFD-subalgebra $N\subset M$, the Hilbert bimodule $_NL^2M_Q$ has a non-zero coarse part. 
If moreover $N$ is quasi-regular in $M$, then $N$ is coarse to $Q$. 
\endproclaim
\noindent
{\it Proof}. This is trivial by Theorem 2.5 and Lemma 2.7.  
\hfill $\square$

\hskip.05in 

For the next result we'ill use the following terminology from [P19a]. An inclusion of factors $\Cal N \subset \Cal M$ 
is ergodic if the action $\Cal U(\Cal N)\curvearrowright^{\text{Ad}}\Cal M$ is ergodic 
(i.e., if $uxu^*=x$ for some $x\in \Cal M$ and all unitary elements $u\in \Cal U(\Cal N)$ then $x$ is necessarily  a scalar multiple of 1). 
The inclusion $\Cal N\subset \Cal M$ is MASA-ergodic if there exists an abelian $^*$-subalgebra $A\subset \Cal N$ that's MASA 
in $\Cal M$, equivalently $A'\cap \Cal M=A$. The inclusion is $R$-ergodic, if there exists a copy of the hyperfinite II$_1$ factor $R\subset \Cal N$ 
that's ergodic in $\Cal M$. And it is AFD-ergodic if there exists an AFD-subalgebra $\Cal R\subset \Cal N$ such that  $\Cal R'\cap \Cal M=\Cal Z(\Cal R)$. 
Also, recall from (Theorem 1.2 in [P19a]) that if an inclusion of factors $\Cal N \subset \Cal M$ is MASA-ergodic, then it is $R$-ergodic. 

\proclaim{2.10. Corollary} Assume $Q\subset M$ is an ergodic compressible inclusion of $\text{\rm II}_1$ factors and denote 
$M_1=\langle M, e_Q \rangle$. Then $M\subset M_1$ is ergodic but not AFD-ergodic, so in particular it is not MASA-ergodic, nor $R$-ergodic. 

More generally, we have: 
\vskip.05in

$1^\circ$ For any Hilbert space $\Cal H_0$, the inclusion $\Cal M:=M\overline{\otimes} \Cal B(\Cal H_0)   
\subset M_1\overline{\otimes} \Cal B(\Cal H_0)$ $=:\Cal M_1$ is ergodic but not AFD-ergodic. 

$2^\circ$ For any  factor $\Cal M_0$, the inclusion $\Cal M:=M\overline{\otimes} \Cal M_0   
\subset M_1\overline{\otimes} \Cal M_0=:\Cal M_1$ is ergodic but not $R$-ergodic. More generally, 
there exists no AFD subfactor $\Cal R\subset \Cal M$ that's ergodic in $\Cal M_1$. 

\vskip.05in 

For example, this is the case if $Q\subset M$ is an ergodic inclusion of $\text{\rm II}_1$ factors such that $M$ contains a diffuse   
tracial W$^*$-algebra $Q_0\subset M$ that's free independent to $Q$. In particular, an  
inclusion of factors of the form $(Q*Q_0) \overline{\otimes} \Cal M_0 \subset \langle Q*Q_0, e_Q\rangle \overline{\otimes} \Cal M_0$ is ergodic 
but not MASA-ergodic, nor $R$-ergodic, and if $\Cal M_0$ is type $\text{\rm I}$ then it is not even AFD-ergodic. 
\endproclaim 
\noindent
{\it Proof}. We first  take $\Cal M_0=\Cal B(\Cal H_0)$. We let $\Cal M=M\overline{\otimes} \Cal B(\Cal H_0)$, $\Cal M_1=\langle M, e_Q\rangle 
\overline{\otimes} \Cal B(\Cal H_0)$ 
be represented on the Hilbert space $\Cal H=L^2M \overline{\otimes}  \Cal H_0$, 
which we also view as a $\Cal M-M$ bimodule. Note that ${M^{op}}'=\Cal M$ and $\Cal M_1={Q^{op}}'$. 

Assume that $\Cal M$ contains an AFD subalgebra $\Cal R\subset \Cal M$ 
such that $\Cal R'\cap \Cal M_1= \Cal Z(\Cal R)$. Note that   
this is equivalent to $\Cal R\vee Q^{op}=\Cal Z(\Cal R)'$.

Since $_{\Cal R}\Cal H_M$ is weakly coarse (because $\Cal R$ is AFD) and $Q\subset \Cal M$ is compressible, 
by  part $1^\circ$ of Theorem 2.4 it follows that there exists a normal conditional expectation  $\Phi: \Cal Z(\Cal R)'=\Cal R \vee Q^{op} \rightarrow \Cal R$ 
satisfying $\Phi(xy^{op})=\tau(y)x$, $\forall x\in \Cal R, y\in Q$. This implies $\Cal R$ has a type I direct summand, so there exists a non-zero 
projection $p\in \Cal R$ such that $p\Cal Rp$ is abelian. Thus, 
the inclusion $\Cal M\subset \Cal M_1$  follows MASA-ergodic, which by [P19a] implies that it is $R$-ergodic. But then taking $\Cal R=R$ in the above,  
it follows that there exists a normal conditional expectation from $\Cal Z(R)'=\Cal B(\Cal H)$ onto $R$. But this is a contradiction, because $\Cal R$ 
is a diffuse (i.e. without atoms) W$^*$-algebra. 

Assume now $\Cal M_0$ is an arbitrary factor and denote 
$\Cal M=M\overline{\otimes} \Cal M_0$, $\Cal M_1=\langle M, e_Q\rangle 
\overline{\otimes} \Cal M_0$.  Consider the Hilbert space $\Cal H=L^2M\overline{\otimes} L^2\Cal M_0=L^2\Cal M$ endowed with its  $\Cal M-\Cal M$ 
bimodule structure. 

Note that if we take 
$B=1\otimes \Cal M_0$ then the compressing unitaries $\Cal U=\Cal U(M)$ for $Q\subset M$ are also 
compressing $\Cal Q=Q\overline{\otimes}\Cal M_0\subset \Cal M$ relative to $B$.  Also, they commute with $B$ and they leave invariant 
the expectation $E=\tau_M \otimes id_{\Cal M_0}$ of $\Cal Q$ onto $B$. 

Assume now that $\Cal R\subset \Cal M$ is an AFD factor such that $\Cal R'\cap \Cal M_1=\Bbb C$. Since $_{\Cal R}\Cal H_{\Cal M}$ is weakly coarse  
(because $\Cal R$ is AFD), all 
the hypothesis in Theorem 2.5 are satisfied. We thus get a normal  conditional expectation $\Phi$ from $\Cal R\vee \Cal Q^{op}$ to $\Cal R \vee B^{op}$. 
Since $\Cal M_1=(\Cal Q^{op})'$, it follows that $(\Cal R\vee \Cal Q^{op})'=\Cal R'\cap \Cal M_1=\Bbb C$ and hence  
$\Cal R\vee \Cal Q^{op}=\Cal B(\Cal H)$.  Thus,      
$\Phi$ is a normal conditional expectation from  $\Cal B(\Cal H)$ onto $\Cal R\vee B^{op}$. This  forces $\Cal R\vee B^{op}$ to be atomic, implying that  
its commutant $(\Cal R\vee B)'\cap \Cal B(\Cal H)$ is atomic as well. Since ${B^{op}}'=(1\otimes \Cal M_0^{op})'=\Cal M_0\overline{\otimes} \Cal B(L^2M)=
\Cal M \vee (M\otimes 1)^{op}$, 
this means that the inclusion $\Cal R\subset \Cal M \subset  \tilde{\Cal M}:=\Cal M \vee (M\otimes 1)^{op}$ has atomic relative commutant  
$\Cal R'\cap \tilde{\Cal M}$. Denoting this atomic algebra  by $\Cal B$, note that we have  $\Cal B \supset M^{op}$, with the commutant of $M^{op}$ in 
$\Cal B$ of type II. But $(M^{op})'\cap (\Cal M \vee M^{op})=\Cal M$, so we have at the same time $\Cal R'\cap \Cal M=\Bbb C$ and $\Cal R'\cap \Cal M
\supset (M^{op})'\cap \Cal B\neq \Bbb C$, contradiction. 
\hfill $\square$

\vskip.05in 

By the free independence theorem in ([P92]), given any separable W$^*$-subalgebras $Q, N$ in an ultraproduct II$_1$ factor $M$, 
there exists a Haar unitary $u\in M$ such that $\{u\}''$ is free independent to both $Q$ and $N$. This implies in particular that $uQu^*\perp N$, 
so $_NL^2M_Q$ has a non-zero coarse part. Thus, any two separable W$^*$-subalgebras of $M$ have a non-zero coarse part. 
Corollary 2.9 above allows us to deduce that if $N$ is AFD, then this is valid for $N$ non-separable as well. It also implies that an AFD $N\subset M$ 
can never be quasi-regular in $M$, 
thus strengthening a result in ([P81b]) showing that ultrapower factors do not have Cartan subalgebras. 

\proclaim{2.11. Corollary} Let $M$ be an ultraproduct $\text{\rm II}_1$ factor. 

$1^\circ$ Given any W$^*$-subalgebras $Q, N\subset M$ with $Q$ separable and $N$ AFD,  the Hilbert $N-Q$ bimodule $_NL^2M_Q$ has a non-zero coarse part. 

$2^\circ$ $M$ does not admit  any diffuse amenable quasi-regular W$^*$-subalgebra. 
\endproclaim
\noindent
{\it Proof}.   $1^\circ$ If $N$ is AFD, then $_NL^2M_M$ is weakly coarse. Since $Q\subset M$ is separable, by Corollary  1.5.4$^\circ$ it is compressible, so Corollary 2.9 applies to get 
that $_NL^2M_Q$ has a non-zero coarse part.    

2$^\circ$ If a W$^*$-subalgebra $N$ of $M$ 
would be AFD and quasi-regular, then by 1.5.4$^\circ$ and 2.9, $_NL^2M_Q$ would follow coarse for any separable 
W$^*$-subalgebra $Q\subset M$. But taking $Q$ to be diffuse included into $N$, gives a contradiction, 
because $_NL^2M_Q$ contains $_NL^2N_Q$, which is not coarse. 
\hfill $\square$



\heading 3. Further remarks and problems
\endheading

\noindent
{\bf 3.1. Compressibility and mean ergodicity of inclusions}. Let $\Cal M$ be a W$^*$-algebra 
and $\Cal T$ be a locally convex topology on $\Cal M$, given by a family of seminorms that are continuous with respect to the operator 
norm $\| \ \|$ on $\Cal M$ (such as the wo, so, or operator norm topology itself).  

If $\Cal G$ is a group 
of automorphisms of $\Cal M$, then we say that $\Cal G \curvearrowright \Cal M$ 
has the $\Cal T$-averaging property if the $\Cal T$-closure of the convex hull of $\{\theta(x) \mid \theta \in \Cal G\}$ has non-empty intersection with the fixed point algebra 
$\Cal M^{\Cal G}:=\{y\in \Cal M \mid \theta(y)=y, \forall \theta \in \Cal G\}$. Note that this property only depends on the set of 
 $\Cal T$-continuous functionals on $\Cal M$, in particular wo-mean ergodicity is same as so-mean ergodicity. 

If $\Cal N \subset \Cal M$ is a W$^*$-subalgebra, then we say that $\Cal N\subset \Cal M$ has the $\Cal T$-averaging property if the Ad-action 
$\Cal U(\Cal N)\curvearrowright^{\text{\rm Ad}} \Cal M$ is $\Cal T$-averaging. 

Note that the $\| \ \|$-averaging of a  W$^*$-inclusion $\Cal N\subset \Cal M$ amounts to what one usually calls 
the relative Dixmier property. It holds true for any ``trivial'' inclusion $\Cal M \subset \Cal M$ by Dixmier's averaging theorem  (see Ch. III, Sec. 5 in [D57]) 
and more generally, by [P96],  for any W$^*$-inclusion $\Cal N\subset \Cal M$ of finite Pimsner-Popa index (as defined in [PP83]). Moreover, if $N\subset M$ is an 
inclusion of II$_1$ factors with $N$ separable, then by (Corollary 4.1 in [P96]), $N\subset M$ has the relative Dixmier property 
(or, in the present terminology, $N\subset M$ is $\| \ \|$-averaging) if and only if  $N$ has finite Jones index in $M$.  
Hence, if $\Cal N\subset \Cal M$ is a  W$^*$-inclusion of finite index, then it is $\Cal T$-averaging for any $\Cal T$ as considered above. 

Note also that the wo-averaging property of a W$^*$-inclusion $\Cal N\subset \Cal M$ has been called ``weak relative Dixmier property'' 
in (Definition 1.1 of [P99]), a terminology that has since then  been adopted  in  several other papers. 
But we should  point out that the idea of ``pushing'' elements $x$ of an ambient algebra $\Cal M$ into the relative commutant 
of a subalgebra $\Cal N \subset \Cal M$ by taking weak limits of averaging of $x$ by unitaries of $\Cal N$ is in fact due to J. Schwartz in ([Sc63]), 
playing a key role in developing the concept of W$^*$-amenability in several equivalent ways. Thus, for a W$^*$-inclusion of the form  $\Cal N \subset \Cal B(\Cal H)$, 
the fact that  $\overline{\text{\rm co}}^w \{uxu^* \mid u\in \Cal U(\Cal N)\}\cap  (\Cal N'\cap \Cal B(\Cal H))\neq \emptyset$, $\forall x\in \Cal B(\Cal H)$,  
amounts to $\Cal N$ satisfying the proprety (P) of Schwartz, later shown equivalent to injectivity (respectively amenability) of $\Cal N$, and also  to 
$\Cal N$ being AFD. 

When referring to an irreducible (or ergodic) 
inclusion of factors $\Cal N\subset \Cal M$, we will say that $\Cal N\subset \Cal M$  is $\Cal T$-mean ergodic 
when it has the $\Cal T$-averaging property. 
This terminology has already been used in [P19a]   
(but in the form ``MV-ergodicity'', and only applied for $\Cal T$ the wo-topology) and is meant to 
emphasize the dynamic aspect of this property that's reminiscent of von Neumann's mean ergodic theorem (see also page 973 in [P19a] for a discussion along 
these lines). 

In this same spirit, the compressibility of a W$^*$-subalgebra $Q\subset \Cal M$ amounts to the action 
$\Cal U(\Cal M)\curvearrowright^{\text{\rm Ad}} \Cal M\overline{\otimes} \Cal B(\ell^2\Bbb N)$ being uniformly $\| \ \|$-averaging 
on the unit ball of $Q \overline{\otimes}\Cal B(\ell^2\Bbb N)$. This property was shown in Corollary 1.5 to hold true whenever 
$\Cal M=M$ is tracial and contains a trace $0$ unitary that's free independent to $Q$. For instance,  
for a free group factor $M=L\Bbb F_n$, $2 \leq n \leq \infty$, the action $\Cal U(M)\curvearrowright^{\text{\rm Ad}} M\overline{\otimes} \Cal B(\ell^2\Bbb N)$ 
is uniformly $\| \ \|$-averaging (or compressible) in all ``abelian directions'' $A_g:=\{u_g\}''$, $g\in \Bbb F_n$. 
Also, by Corollary 1.7 we have that any interpolated free group factor $M=L\Bbb F_t, 1< t \leq \infty$, admits an increasing sequence of 
irreducible subfactors $M_m$ exhausting $M$, such that $M$ is uniformly $\| \ \|$-averaging (or compressible) on each $M_m\overline{\otimes} \Cal B(\ell^2\Bbb N)$. 

A related property for a II$_1$ factor $M$ is the $\| \ \|$-averaging of $\Cal U(M)\curvearrowright M\overline{\otimes} \Cal B(\ell^2\Bbb N)$ 
(in other words, the relative Dixmier property for the inclusion $M\subset M\overline{\otimes} \Cal B(\ell^2\Bbb N)$).  
We will call this  the $\| \ \|_{cb}$-averaging property of $M$. 
Note that despite results in [P22] showing that in any II$_1$ factor $M$ one can simultaneously ``$\varepsilon$-push'' to the scalars arbitrarily large finite sets of elements in  
$(M)_1$ by averaging 
with  $n=n(\varepsilon)$ unitaries in $M$, this does not  seem to entail  $M$ is $\| \ \|_{cb}$-averaging. 

It would be interesting 
to know whether the hyperfinite II$_1$ factor and the free group factors are $\| \ \|_{cb}$-averaging.   
Let us however note  that the ultrapower of any II$_1$ factor does have the property. More generally, we have:

\proclaim{3.1.1. Proposition}  Any ultraproduct $\text{\rm II}_1$ factor $\Cal M=\underset{n \rightarrow \omega}\to{\Pi} M_n$, 
with $M_n$ a sequence of finite factors satisfying $\text{\rm dim}(M_n)\nearrow \infty$ and $\omega$ a free ultrafilter on $\Bbb N$, is $\| \ \|_{cb}$-averaging.  
\endproclaim
\noindent
{\it Proof}. We have to prove that given any $x=(x_{kl})_{k, l\in \Bbb N}\in \Cal M\overline{\otimes} \Cal B(\ell^2\Bbb N)$ and any $\varepsilon>0$, 
there exists a finite set $\Cal U_0\subset \Cal U(M)$ such that 
$$
\|\frac{1}{|\Cal U_0|} \sum_{u\in \Cal U_0} uxu^*-E_{1\otimes \Cal B(\ell^2\Bbb N)}(x)\|\leq \varepsilon. \tag 1
$$ 

Since $\Cal X=\{x_{kl} \mid k,l \in \Bbb N\}$ is countable, by [P92] there exists a Haar unitary $v\in \Cal M$ such that $A=\{v\}''$ is free independent to $\Cal X''\subset \Cal M$. 
This implies that for any $K\geq 1$, the truncation $x^K=(x_{kl})_{1\leq k,l\leq K}\in \Cal M\otimes\Bbb M_K(\Bbb C)$ generates a $^*$-algebra that's 
free independent to $A$ relative to $1\otimes \Bbb M_K(\Bbb C)$. By Proposition 1.4, it follows that if $u\in A$ is a Haar unitary 
and we let $\Cal U_0=\{1, u, ..., u^{n-1}\}$ for some $n$ satisfying $2\sqrt{n-1}/n$, then 
$$
\|\frac{1}{|\Cal U_0|} \sum_{u\in \Cal U_0} ux^Ku^*-E_{1\otimes \Bbb M_K(\Bbb C)}(x^K)\|\leq \varepsilon. \tag 2
$$
Since the choice of $\Cal U_0$ is independent of $K$, taking $K\rightarrow \infty$ in $(2)$  implies that $(1)$ holds true as well. 
\hfill $\square$

\vskip.1in

\noindent
{\bf 3.2. The wFC property and compressibility}. The free complementation (FC) problem asks whether any maximal amenable W$^*$-subalgebra 
$Q$ of a free group factor $M=L\Bbb F_n$, $2\leq n \leq \infty$, admits a ``free complement'' in $M$, i.e., there exists a W$^*$-subalgebra 
$N\subset M$ such that $M=Q\vee N\simeq Q*N$. We refer to (page 3106 in [P18]; Remark 1.4.2 in [BP23];  Section 4.1 in [BDIP23]) 
for more detailed discussions on this problem (see also [BDH24] for some 
recent progress in this direction). The following weaker property, that we call wFC, has been conjectured to hold true:  

\vskip.05in 
\noindent
3.2.1. {\it  wFC conjecture}. Given any amenable W$^*$-subalgebra $Q$ in a free group factor $M$, there exists a diffuse 
abelian W$^*$-subalgebra $A\subset M$ that's free independent to $Q$. 

\vskip.05in

Note that if $Q\subset M$ is diffuse then the existence of a diffuse W$^*$-subalgebra $A\subset M$ that's free independent to $Q$ 
is equivalent to the existence of a trace zero selfadjoint unitary in $M$ that's free independent to $Q$. Given a 
W$^*$-subalgebra $B$ in a II$_1$ factor $M$, let us denote $F_0(B)=\{u=u^*\in \Cal U(M)\mid \tau(u)=0, \{u\}$ free independent to $B\}$. 
Thus, the wFC property amounts to $F_0(Q)\neq \emptyset$, for any amenable $Q\subset M$. The following strengthening of this condition 
should in fact hold true: 

\vskip.05in 
\noindent
3.2.2. {\it Strengthened wFC conjecture}. Given any maximal amenable W$^*$-subalgebra $Q$ in a free group factor $M$, 
the weak closure of the linear span of $F_0(Q)$ is equal to $M\ominus Q$.

\vskip.05in

It is easy to see that if a free group factor satisfies 3.2.2 then it satisfies the coarseness conjecture ([H15], [P18]), 
stating that any maximal amenable subalgebra $Q\subset M$ is coarse. In fact, for $M$ to satisfy the coarseness conjecture,   it is sufficient 
that $\overline{\text{\rm sp} QF_h(Q)Q}^w =M\ominus Q$, where $F_h(Q):=\{h\in M_h \mid \tau(h)=0, \{h\}$ free independent to $Q\}$. 
The coarseness conjecture is known to imply 
the Peterson-Thom conjecture (see Proposition 5.3 in [P18]) and both were settled in the affirmative in ([H20], [BeCa22], [BoCo23]).  Thus, in addition to 
providing a striking structural property of free group factors,  3.2.2 above could lead to a new approach to (and a strenghtening of) these results.

On the other hand,  by  Corollary 1.5, if $M$ satisfies 3.2.1, then any amenable W$^*$-subalgebra in $M$ follows compressible. So 
in particular,  the following weaker property would hold true as well: 

\vskip.05in 
\noindent
3.2.3. {\it Compressibility conjecture}. Any amenable W$^*$-subalgebra $Q$ of a free group factor $M$ is compressible. Equivalently, 
the action $\Cal U(M)\curvearrowright^{\text{\rm Ad}} M\overline{\otimes} \Cal B(\ell^2\Bbb N)$ is uniformly $\| \ \|$-averaging 
on $(Q\overline{\otimes} \Cal B(\ell^2\Bbb N))_1$ for any amenable W$^*$-subalgebra $Q$ of $M$.

\vskip.1in

\noindent
{\bf 3.3. Condensation of freeness phenomena}. There is an approach to the wFC problem which at an initial stage seems to suggest that any separable 
II$_1$ factor $M$ (so in particular any free group factor) contains a diffuse abelian W$^*$-subalgebra  $A\subset M$ that admits no 
elements free to it, $F_0(A)=\emptyset$.  This is based on an ``AFD-percolation approach'' to the problem, which 
in this specific case consists in constructing $A$ recursively as 
an inductive limit of finite partitions, $A_n\nearrow A$, such that at each step $n$ the finite dimensional algebra $A_n$  
is chosen so that ``more and more'' of the projections in $\Cal P_0(M):=\{p\in \Cal P(M)\mid \tau(p)=1/2\}$ are being prevented from being free from $A_n$. 
It is  in fact sufficient to do this with only a dense subset of projections  $\{p_m\}_m\subset \Cal P_0(M)$, provided one can ``spoil'' freeness 
in some uniform way. An example of a  local property that allows such recursive machinery to work is the following:  

\vskip.05in 
\noindent
3.3.1. {\it $\exists c_0>0$ such that $\forall B_0\subset M$ finite partition, $\forall p_0\in F_0(B_0)$, $\exists B_0\subset B_1\subset M$ finite partition 
with the property that $\|p_1-p_0\|_2\geq c_0$, $\forall p_1\in F_0(B_1)$}.

\vskip.05in 

Assume $M$ satisfies the above condition. Let $\{p_m\}_m \subset \Cal P_0(M)$  be a sequence of projections that's dense in $\Cal P_0(M)$ 
and $\{x_n\}_n\subset (M)_1$ be $\| \ \|_2$-dense in $(M)_1$. 
We construct recursively an increasing sequence of finite partitions $A_n\subset M$ such that at each step $n$, the $\| \ \|_2$-distance 
between $p_n$ and $F_0(A_n)$, $\delta(p_n, F_0(A_n))$, satisfies:
$$
\delta(p_n, F_0(A_n))\geq c_0,  \tag 1
$$
and also 
$$
\|E_{A_n}(x_k)-E_{A_n'\cap M}(x_k)\|_2\leq 2^{-n}, 1\leq k \leq n. \tag 2
$$

The role of condition $(1)$ is to ``spoil'' freeness between $p_n$ and $A_n$ (and thus all ``future'' $A_m\supset A_n$ as well!), while the role of $(2)$ is to make $A_n$  
more and more a MASA in $M$. 

If we reached step $n$, then $A_{n+1}$ is constructed by first applying $(1)$ to $B_0=A_n$ to get a refinement  $B_1=A^{0}_{n+1}$ 
such that $\delta(p_{n+1}, F_0(A_{n+1}^0))\geq c_0$. Then by using (Lemma 1.2 in [P81a]) one gets a refinement $A_{n+1}\supset A^{0}_{n+1}$ 
such that $\|E_{A_{n+1}}(x_k)-E_{A_{n+1}'\cap M}(x_k)\|_2\leq 2^{-n-1}$, $1\leq k \leq n+1$. 

Finally, one takes $A$ to be the weak closure of $\cup_n A_n$. By $(1)$, if $F_0(A)$ would be non-empty, we would have 
$$
\delta(p_n, F_0(A))\geq \delta(p_n, F_0(A_n))\geq c_0, \forall n, \tag 3
$$ 
which by the density of $\{p_n\}_n$ in $\Cal P_0(M)$ implies $\delta(p, F_0(A))\geq c_0$, $\forall p\in \Cal P_0(M)$, a contradiction. Thus, 
$F_0(A)=\emptyset$. At the same time, by (1.2 in [P81a]) condition $(2)$ insures that $A=A'\cap M$, i.e., $A$ is a MASA in $M$. 

We have thus shown:

\proclaim{3.3.2. Proposition} If a separable  $\text{\rm II}_1$ factor $M$ satisfies condition $3.3.1$ then $M$ admits a MASA $A\subset M$ 
with the property that $F_0(M)=\emptyset$. In particular, if a free group factor satisfies $3.3.1$ then the wFC conjecture $($and thus also the FC conjecture$)$ fails. 
\endproclaim

Condition $3.3.1$  seems in fact too strong for the purpose of ``spoiling'' freeness.
 Ideally, such a condition should be so that $M$ satisfies it if and only if $M$ is non-wFC.  
The negation of 3.3.1 though is the following 

\vskip.05in 
\noindent
3.3.1'. {\it $\forall c>0$, $\exists B_0\subset M$ finite partition, $\exists p_0\in F_0(B_0)$, 
such that given any refinement  $B_0\subset B_1\subset M$ of $B_0$ one has $\delta(p_0, F_0(B_1))< c$}. 
\vskip.05in 
\noindent
and this condition does not seem to imply in a straightforward way that $M$ satisfies the wFC conjecture. Indeed, it is not clear at all 
that 3.3.1' is sufficient to insure the possibility of choosing   
$p_n\in F_0(A_n)$ in a Cauchy-manner, for some $A_n \nearrow A$ (a fact that would provide a projection $p\in F_0(A)$ 
by simply taking the limit of the $\{p_n\}_n$). What does allow such a choice is the following strengthening of 3.3.1': 

\vskip.05in 
\noindent
3.3.3. {\it $\forall c>0$, $\exists \delta>0$ 
such that given any $B_0\subset M$ finite partition with  mesh $\delta(B_0)$ less than $\delta $, one has $\delta(F_0(B_0), F_0(B_1))< c$, 
for any finite partition $B_1\subset M$ refining $B_0$.} 
\vskip.05in 

Indeed, it is easy to see that if a II$_1$ factor $M$ satisfies 3.3.3 then it satisfies wFC. To show this, one constructs recursively 
an increasing sequence of partitions $A_n$, with $\delta(A_n)\leq 2^{-n}$, together with a choice of some $p_n\in F_0(A_n)$ 
such that $\|p_{n-1}-p_n\|_2\leq 2^{-n}$, at each step $n\geq 1$. This gives rise to a Cauchy sequence $\{p_n\}_n$ whose  limit $p$  
lies in $F_0(A)$, where $A=\overline{\cup_n A_n}^w$. 
We will call 3.3.3 the {\it condensation of freeness} property. 
So the following problem seems of interest: 

\vskip.05in 
\noindent
3.3.4. {\it Condensation of freeness problem}. Do the free group factors (more generally the interpolated 
free group factors)  satisfy the condensation of freeness property ?

\vskip.1in

\noindent
{\bf 3.4. A motivation related to Connes' embedding problem}. We initiated the present work some time ago, 
as an approach to Connes embedding (CE) problem. 

The starting point  was the following intuition. Assuming $M$ is a II$_1$ factor 
with an irreducible subfactor $Q$ for which the action $\Cal U(M)\curvearrowright^{\text{\rm Ad}}M$ is uniformly 
ergodic on $(Q)_1$ in some appropriate topology (a property which 
we'll call ``compressibility'') and such that there exists a hyperfinite subfactor $N\subset M$ 
with the property that $N\vee  Q^{op}=\Cal B(L^2M)$, then by Kaplanski's density theorem one can 
so-approximate any finite set of elements  in $(M)_1$ by elements of the form $T=\sum_j x_j \cdot y^{op}_j\in (N \vee_{Alg} Q^{op})_1$. 
Evaluating such operators $T$ at $\langle \ \cdot \ \hat{u}_k, \hat{u}_k\rangle$, $1\leq k \leq n$, where $u_j\in \Cal U(M)$ are chosen so that 
$\frac{1}{n}\sum_j u_jyu_j^*$ is close to scalars uniformly for all $y\in (Q)_1$ (by ``compressibility''), one gets that the so-approximations $T$ are in fact 
wo-close to being in $N=N \vee 1$, thus giving some weak approximation of the elements in $M$ by  elements in the AFD algebra $N$, 
and hence by matrix algebras. 

Such  finite dimensional wo-approximation though is immediately seen not to be possible, or else $M$ itself would follow amenable (as the wo-approximation happens within $M$). 

But one can amplify $M$ by $\Cal B(\ell^2\Bbb N)$ acting on the Hilbert space $\Cal H=L^2M\otimes \ell^2\Bbb N$, viewed as a left module 
over $\Cal M=M\overline{\otimes} \Cal B(\ell^2\Bbb N)$ and right module over $M$. If there exists $Q\subset M$  that's ``compressible'' 
and if there exists $N\simeq R$ in $\Cal M$ such that $N\vee Q^{op}=\Cal B(\Cal H)$ 
then the same argument would give some weak approximations of $M$ by matrix subalgebras of $N$ that are this time ``exterior'' 
to $M$ (as $N$ lies in the larger algebra $\Cal M$). 

Theorem 2.5 and its Corollary 2.9  show that even this is in fact not possible, at least not with the notion of ``compressibility'' as defined in 2.1: 
such an assumption on $Q$ makes it impossible for  an $R$-tight complement of $Q$ to exist in $\Cal M=M\overline{\otimes} \Cal B(\ell^2\Bbb N)$. Thus, this approach   
doesn't seem to be usable towards obtaining  criteria  for a II$_1$ factor to be Connes embeddable.

\vskip.1in 
\noindent
{\bf  Conflict of Interest and Data Statement.} The author certifies that he has no conflict interests of any kind related to 
the content of this manuscript. 
Data sharing is not applicable to this article as no new data were created or analyzed in this work.

\head  References \endhead

\item{[AkO76]} C.A. Akemann, P.A. Ostrand: {\it Computing norms in group} C$^*$-{\it algebras}. Amer. J. Math. {\bf 98} (1976), 1015-1047

\item{[A95]} C. Anantharaman: {\it Amenable correspondences and approximation properties 
for von Neumann algebras}, Pacific J. Math., {\bf 171} (1995), 309-343. 

\item{[AP17]} C. Anantharaman, S. Popa: ``An introduction to II$_1$ factors'', \newline www.math.ucla.edu/$\sim$popa/Books/IIun-v17.pdf

\item{[BeCa22]} S. Belinschi, M. Capitaine: {\it Strong convergence of tensor products of independent G.U.E. matrices}, arXiv:2205.07695 

\item{[BoCo23]}  C. Bordenave, B. Collins:  {\it Norm of matrix-valued polynomials in random unitaries and permutations}, arXiv:2304.05714

\item{[BDH24]} N. Boschert, E. Davis, P. Hiatt: {\it A Class of Freely Complemented von Neumann Subalgebras of $L\Bbb F_n$}, to appear in 
Intern. Math. Res. Notices, arXiv:2411.05136

\item{[BDIP23]} R. Boutonnet, D. Drimbe, A. Ioana, S. Popa: {\it Non-isomorphism of $A^{*n}, 2\leq n \leq \infty$, 
for a non-separable abelian von Neumann algebra $A$}, 
GAFA 2024, \newline  https://doi.org/10.1007/s00039-024-00669-8, 2024  (arXiv:2308.05671). 

\item{[BP23]} R. Boutonnet, S. Popa: {\it Maximal amenable MASAs of radial type in the free group factors},  
to appear in Proc. AMS,  arXiv:2302.13355

\item{[BrO08]} N. Brown, N. Ozawa: ``C$^*$-algebras and finite dimensional approximation'', Grad. Studies in Math. Vol. {\bf 88}, 
AMS  Providence, Rhode Island, 2008.

\item{[C76]} A. Connes: {\it Classification of injective factors}, Ann. of Math., {\bf 104} (1976), 73-115.

\item{[DiP24]} C. Ding, J. Peterson: {\it Biexact von Neumann algebras}, arXiv:2309.10161

\item{[D57]} J. Dixmier: ``Les alg\' ebres d'op\' erateurs sur l'espace Hilbertien (Alg\' ebres de von Neumann)'', Gauthier-Villars, Paris 1957.

\item{[Dy93]} K. Dykema: {\it Free products of hyperfinite von Neumann algebras and free dimension}, Duke Math. J. {\bf 69} (1993), 97-119. 

\item{[Dy94]}  K. Dykema: {\it Interpolated free group factors}, Pacific J. Math. {\bf 163} (1994),123-135. 

\item{[EL77]} E. Effros, C. Lance: {\it Tensor products of operator algebras}, Advances in Math., {\bf 25} (1977), 1-34. 

\item{[GP96]} L. Ge, S. Popa: {\it On some decomposition properties
for factors of type} II$_1$, Duke Math. J., {\bf 94} (1998), 79-101.

\item{[Ha15]} B. Hayes: 1-{\it Bounded entropy and regularity problems in von Neumann algebras}. International Math. Res. Notices, {\bf 1} 
(2018), 57-137 (arXiv:1505.06682)

\item{[Ha20]} B. Hayes: {\it A random matrix approach to the Peterson-Thom conjecture}. Indiana University Mathematics Journal. 2022, {\bf 71} (2022), 
1243-1297 (arXiv:2008.12287) 

\item{[HoI23]} C. Houdayer, A. Ioana: {\it Asymptotic freeness in tracial ultraproducts},  Forum Math. Sigma {\bf 12} (2024), Paper No. e88, 22 pp, arXiv:2309.15029

\item{[IT23]} T. Hui, A. Ioana: {\it On existentially closed} II$_1$ {\it factor}, J. Funct. Anal. {\bf 286} (2024), Paper No. 110264 (arxiv.org/pdf/2306.00474)

\item{[Ke59]} H. Kesten: {\it Symmetric random walks on groups}, Trans. Amer. Math. Soc. {\bf 92} (1959), 336-354. 

\item{[M23]} A. Marrakchi: {\it Kadison's problem for type} III {\it subfactors and the bicentralizer conjecture}, 
Inventiones Math. {\bf 239} (2025), 79-163

\item{[PiP83]} M. Pimsner, S. Popa: {\it Entropy and index for subfactors}, Ann. Sci. Ec. Norm. Sup. {\bf } (1986), 57-106. 

\item{[P81]} S. Popa:  {\it On a problem of R.V. Kadison on maximal abelian *-subalgebras in factors}, Invent. Math., {\bf 65} (1981), 269-281.

\item{[P81b]} S. Popa: {\it Orthogonal pairs of *-subalgebras in
finite von Neumann algebras}, J. Operator Theory, {\bf 9} (1983),
253-268.

\item{[P82]} S. Popa: {\it Maximal injective subalgebras in factors associated with free groups}, Advances in Math., {\bf 50} (1983), 27-48.

\item{[P86]} S. Popa: {\it Correspondences}, INCREST preprint 56/1986, unpublished, see \newline https://www.math.ucla.edu/$\sim$popa/preprints.html 

\item{[P90]} S. Popa: {\it Markov traces on Universal Jones algebras
and subfactors of finite index}, Invent.  Math., {\bf 111} (1993),
375-405.

\item{[P92]} S. Popa: {\it Free independent sequences in type} II$_1$ 
{\it factors and related problems}, Asterisque {\bf 232} (1995), 187-202.

\item{[P96]} S. Popa: {\it The relative
Dixmier property for inclusions of von Neumann
algebras of finite index}, Ann. Sci. Ec. Norm. Sup.,
{\bf 32} (1999), 743-767.

\item{[P01]} S. Popa: {\it On a class of type} II$_1$ {\it factors with
Betti numbers invariants}, Ann. of Math {\bf 163} (2006), 809-899
(MSRI preprint 2001-024; math.OA/0209310).  

\item{[P03]} S. Popa: {\it Strong Rigidity of} II$_1$ {\it Factors Arising from Malleable Actions of $w$-Rigid Groups} I, Invent. Math.,
{\bf 165} (2006), 369-408. (math.OA/0305306).

\item{[P13]} S. Popa: {\it Independence properties in subalgebras of ultraproduct} II$_1$ {\it factors}, Journal of Functional Analysis 
{\bf 266} (2014), 5818-5846 (math.OA/1308.3982) 

\item{[P18]} S. Popa: {\it Coarse decomposition of} II$_1$ {\it factors},  Duke Math. J. {\bf 170} (2021) 3073 - 3110,  math.OA/1811.09213

\item{[P19a]} S. Popa: {\it On ergodic embeddings of factors}, Communication in Math. Physics,  {\bf 384} (2021), 971-996 (math.OA/1910.06923)

\item{[P19b]} S. Popa:  {\it Tight decomposition of factors and the single generation problem}, Journal of Operator Theory, {\bf 85} (2021), 277-301 (arXiv:1910.14653)

\item{[P22]} S. Popa: {\it On the paving size of a subfactor}, Pure and Applied Math Quarterly (volume dedicated to Vaughan Jones), 
{\bf 19} (2023), 2525-2536 (arXiv:2210.04396)

\item{[PS18]} S. Popa, D. Shlyakhtenko:  {\it Representing the interpolated free group factors as group factors},  
Groups, Geometry, and Dynamics, {\bf 14} (2020), 837-855  \newline 
(math.OA/1805.10707)

\item{[PV11]} S. Popa, S. Vaes: {\it Unique Cartan decomposition for} II$_1$ 
{\it factors arising from arbitrary actions of free groups}, Acta Mathematica, {\bf 194} (2014), 237-284 \newline (math.OA/1111.6951)

\item{[PV14]}  S. Popa, S. Vaes: {\it On  the optimal paving over MASAs in von Neumann algebras},  Contemporary Mathematics
Volume {\bf 671}, R. Doran and E. Park editors, American Math Society, 2016, pp 199-208. 

\item{[R94]} F. Radulescu: {\it Random matrices, amalgamated free products and subfactors of the von Neumann algebra of a free group, of noninteger index}, 
Invent. Math. {\bf 115} (1994), 347-389. 

\item{[S71]} S. Sakai: ``C$^*$-algebras and W$^*$-algebras'', Springer-Verlag, 1971. 

\item{[Sc63]} J. Schwartz: {\it Two finite, non-hyperfinite, non-isomorphic factors}, Comm. Pure App. Math.  {\bf 16} (1963) 19-26. 

\item{[T79]} M. Takesaki: ``Theory of operator algebras'', I, Springer-Verlar, New York 1979.
 
\item{[Vo88]} D. Voiculescu: {\it Circular and semicircular systems and free product factors}, Prog. in Math. {\bf 92}, Birkhauser, Boston, 1990, pp. 45-60.  

\enddocument